\documentclass[twocolumn,amsthm]{autart}  

\usepackage{xcolor}
\usepackage{hyperref}
\usepackage[ocgcolorlinks]{ocgx2}[2019/01/02] 

\usepackage[switch,pagewise,displaymath,mathlines]{lineno}
\modulolinenumbers[1]

\usepackage{xparse}
\usepackage{amsmath,amssymb,amsfonts,amsxtra,bm,mathtools}
\usepackage{csquotes}
\usepackage{ifthen}
\usepackage{tikz}
\usetikzlibrary{calc,positioning,intersections,arrows.meta,matrix,decorations.pathmorphing,fit}
\usepackage{pgfplots}
\pgfplotsset{compat=1.14}

\newcommand{\ie}{i.\,e.,}
\newcommand{\eg}{e.\,g.}
\newcommand{\wrt}{w.\,r.\,t.}
\newcommand{\xt}{\tilde{x}}
\newcommand{\At}{\widetilde{A}}

\renewcommand{\d}{\mathrm{d}}
\newcommand{\pxi}{\bm{\xi}}
\newcommand{\peta}{\bm{\eta}}
\newcommand{\s}{\bm{s}}
\newcommand{\pG}{\bm{G}}
\newcommand{\pH}{\bm{H}}
\newcommand{\pJ}{\bm{J}}
\newcommand{\oz}{\bar{z}}
\newcommand{\K}{\bm{K}}
\newcommand{\J}{\bm{J}}
\newcommand{\as}{\bar{\bm{a}}_{\Sigma}}
\newcommand{\ad}{\bar{\bm{a}}_{\Delta}}
\newcommand{\F}{\bm{F}}
\newcommand{\oxi}{\bar{\xi}}
\newcommand{\oeta}{\bar{\eta}}
\newcommand{\intl}{\int\limits}
\newcommand{\mA}{\mathcal{A}}
\newcommand{\pX}{\bm{X}}

\newcommand{\y}{\check{z}}

\newcommand{\wt}{\tilde{w}}
\newcommand{\yn}{\y_0}
\newcommand{\ytn}{\symb{ytn}}

\newcommand{\xf}{\bar{x}}
\newcommand{\xfr}{\bar{x}^{r}}

\newcommand{\Af}{\check{A}}
\newcommand{\T}{\mathcal{T}}
\newcommand{\Q}{\bm{Q}}
\renewcommand{\P}{\bm{P}}
\newcommand{\Ab}{\bar{A}}

\newcommand{\li}{\lambda_i^{r}}
\newcommand{\ljLow}{
	\lambda_
	{\raisebox{0pt}[\heightof{$\scriptstyle j$}][0pt]{$\scriptstyle j$}}
	^{r}}

\newcommand{\comment}[1]
{}

\newcommand{\mrk}[2][1]
{%
	\ifthenelse{1>2}{%
		{\color{red}{#2}}%
	}{#2}%
}%

\newcommand{\cross}[1]
{}

\newcommand{\utriag}[1][\scriptstyle]{
\tikz[line join=round]{
	\node(a)[inner sep=0pt,outer sep=0pt]{\phantom{$#1 2$}};
	\draw(a.north west)--(a.south east)--(a.north east)--cycle;
}
}
\newcommand{\ltriag}[1][\scriptstyle]{
\tikz[line join=round]{
	\node(a)[inner sep=0pt,outer sep=0pt]{\phantom{$#1 2$}};
	\draw(a.north west)--(a.south west)--(a.south east)--cycle;
}
}

\NewDocumentCommand{\xl}{e{^_}}{
	\IfNoValueTF{#1}
		{
		\IfNoValueTF{#2}
			{x^{l}} 
			{x^{l}_{#2}}
		} 
		{
		\IfNoValueTF{#2}
			{x^{l#1}} 
			{x^{l#1}_{#2}} 
		}
}
\NewDocumentCommand{\xr}{e{^}}{
	\IfNoValueTF{#1}
	{x^{r}}
	{x^{r#1}}
}

\NewDocumentCommand{\xtl}{e{^_}}{
	\IfNoValueTF{#1}
		{
		\IfNoValueTF{#2}
			{\xt^{l}} 
			{\xt^{l}_{#2}}
		} 
		{
		\IfNoValueTF{#2}
			{\xt^{l#1}} 
			{\xt^{l#1}_{#2}} 
		}
}

\NewDocumentCommand{\xtr}{e{^}}{
	\IfNoValueTF{#1}
	{\xt^{r}}
	{\xt^{r#1}}
}

\makeatletter
\newcommand\getheightofnode[2]{%
	\pgfextracty{#1}{\pgfpointanchor{#2}{north}}%
	\pgfextracty{\pgf@xa}{\pgfpointanchor{#2}{south}}
	\addtolength{#1}{-\pgf@xa}%
	\global#1=#1
}

\newcommand\getwidthofnode[2]{%
	\pgfextractx{#1}{\pgfpointanchor{#2}{east}}%
	\pgfextractx{\pgf@xa}{\pgfpointanchor{#2}{west}}
	\addtolength{#1}{-\pgf@xa}%
	\global#1=#1
}
\makeatother

\colorlet{greenCol}{green!70!black}
\colorlet{blueCol}{blue!70!gray}
\colorlet{redCol}{red!70!black!80!white}
\colorlet{fillArea}{gray!10}
\colorlet{colBC2}{greenCol}
\colorlet{colBC1}{blueCol}
\colorlet{colBC3}{redCol}
\tikzset{insn/.style={inner sep=0pt}}

\allowdisplaybreaks[1]

\newcommand{\rightcond}[1]{\scriptsize\setlength{\arraycolsep}{0.5pt}\renewcommand{\arraystretch}{0.8}\begin{array}[t]{rl}
#1
\end{array}}

\usepackage{arrayjobx}[2010/05/03]
\newarray\constNames 
\newarray\constNamesHigh 
\newarray\constNamesK 
\newarray\constNamesSymb 
\newcounter{cnstcnt} 
\newcounter{cnstcntHigh} 
\newcounter{cnstcntK} 
\newcounter{cnstcntSymb} 
\makeatletter

\newcommand{\cDef}[2][\gamma]{%
	\refstepcounter{cnstcnt}
	\protected@write\@auxout{}
		{\string\constNames(\thecnstcnt)={#1}}
	\protected@write\@auxout{}
	{\string\newlabel{#2}{{\thecnstcnt}{\thepage}}}%
	\hypertarget{#2}{#1_{\thecnstcnt}}%
}
\newcommand{\cHDef}[2][\bm{c}]{%
	\refstepcounter{cnstcntHigh}
	\protected@write\@auxout{}
		{\string\constNamesHigh(\thecnstcntHigh)={#1}}
	\protected@write\@auxout{}
	{\string\newlabel{#2}{{\thecnstcntHigh}{\thepage}}}%
	\hypertarget{#2}{#1^{\thecnstcntHigh}}%
}
\newcommand{\kDef}[2][k]{%
	\refstepcounter{cnstcntK}
	\protected@write\@auxout{}
		{\string\constNamesK(\thecnstcntK)={#1}}
	\protected@write\@auxout{}
	{\string\newlabel{#2}{{\thecnstcntK}{\thepage}}}%
	\hypertarget{#2}{#1_{\thecnstcntK}}%
}

\newcommand{\symbDef}[2]{
	\refstepcounter{cnstcntSymb}
	\protected@write\@auxout{}
		{\string\constNamesSymb(\thecnstcntSymb)={#1}}
	\protected@write\@auxout{}
	{\string\newlabel{#2}{{\thecnstcntSymb}{\thepage}}}%
	\hypertarget{#2}{#1}%
}
\newcounter{currConst}
\newcounter{currConstK}
\NewDocumentCommand{\symb}{o m}{%
	\setcounterref{currConst}{#2}
	\IfNoValueTF{#1}{
	\hyperlink{#2}{\constNamesSymb(\value{currConst})}
	}
	{
	\hyperlink{#2}{\constNamesSymb(\value{currConst})_{#1}}
	}
}%

\NewDocumentCommand{\cH}{o m}{%
	\setcounterref{currConst}{#2}
	\IfNoValueTF{#1}{
		\hyperlink{#2}{\constNamesHigh(\value{currConst})^{\thecurrConst}}
	}{
		\hyperlink{#2}{\constNamesHigh(\value{currConst})^{\thecurrConst}_{#1}}
	}
}

\let\oldequation\equation
\let\oldsubequations\subequations
\let\oldalign\align
\let\oldmultline\multline

\let\oldendequation\endequation
\let\oldendsubequations\endsubequations
\let\oldendalign\endalign
\let\oldendmultline\endmultline

\renewenvironment{equation}
  {\linenomathNonumbers\oldequation}
  {\oldendequation\endlinenomath}
\renewenvironment{align}
	{\linenomathNonumbers\oldalign}
	{\oldendalign\endlinenomath}
\renewenvironment{subequations}
	{\linenomathNonumbers\oldsubequations}
	{\oldendsubequations\endlinenomath}
\renewenvironment{multline}
	{\linenomathNonumbers\oldmultline}
	{\oldendmultline\endlinenomath}

\begin{document}
\begin{frontmatter}

\title{Bilateral backstepping control of coupled \\ linear parabolic PDEs with spatially varying coefficients}


\author[Simon]{Simon Kerschbaum}\ead{simon.kerschbaum@fau.de}~and~
\author[Joachim]{Joachim Deutscher}\ead{joachim.deutscher@uni-ulm.de}
\address[Simon]{Lehrstuhl f\"ur Regelungstechnik, Universit\"at Erlangen-N\"urnberg, 
Cauerstra{\ss}e 7, D--91058 Erlangen, Germany}
\address[Joachim]{Institut für Mess-, Regel- und Mikrotechnik, Universität Ulm, Albert-Einstein-Allee 41, D-89081 Ulm, Germany}  


\begin{abstract}
This paper considers the backstepping state feedback control of coupled linear parabolic PDEs with spatially varying coefficients and bilateral actuation. 
By making use of the folding technique, a system representation with unilateral actuation is obtained, allowing to apply the standard backstepping transformation. 
To ensure the regularity of the solution, the folded system is subject to unusual folding boundary conditions, 
which lead to additional boundary couplings between the PDEs. Therefore, the solution of the corresponding kernel equations determining the transformations is a very challenging problem.
A systematic approach to derive the corresponding integral equations is proposed, allowing to solve \mrk{them} with the method of successive approximations.
By making use of a Volterra and a Volterra-Fredholm transformation, the closed-loop system is mapped into a cascade of stable parabolic systems. This allows a simple proof of 
exponential stability in the $L_2$-norm with the decay rate as design parameter.
The bilateral state feedback stabilization of an unstable system of two coupled parabolic PDEs and the comparison to the application of an unilateral controller demonstrates the results of the paper.
\end{abstract}

\end{frontmatter}

\setlength{\parskip}{0pt}
\setlength{\parindent}{1em}


\section{Introduction}
\subsection{Background and motivation}
In the last decades, the \emph{backstepping method} has emerged as a powerful tool for the boundary control of distributed parameter systems (for an overview, see \cite{Kr08,Vaz17}).
In the scope of parabolic systems, current research is focused on the control of coupled \cross{parabolic} PDEs, where results exist for both constant \cite{Ba15a,Orl17} and spatially varying coefficients \cite{Vaz16a,Deu17} as well as for space and time dependent parameters \cite{Ker19}.
This system class is of great interest for applications, modelling technological processes in chemical and biochemical engineering (see \cite{At74,Ja14}) and crystallisation processes (see \cite{Cr02}).
\mrk{Most existing results for the control of parabolic systems consider the \emph{unilateral} case, which means that the inputs of the system are all located at one boundary.}

In practice, however, there may be the possibility to place actuators at all boundaries of the domain, achieving \emph{bilateral} actuation. Using this additional degree of freedom in the controller design, it is to be expected that the control effort at one input for achieving comparable results as in the unilateral case, can be significantly reduced.
Furthermore, the additional actuation has shown to be useful in the deployment of agents, when modelling multi-agent
systems by continuum models (see \cite{Fre20}). In addition, bilateral control of traffic flow models is considered in \cite{Bek19}.

In \cite{Vaz16c}, the backstepping method was used for the controller and observer design of parabolic PDEs \mrk{on balls of arbitrary dimension} with actuation on the whole boundary. Since in one dimension, the ball actually is an interval, this included bilateral control of 1D parabolic systems. Therefore, \cite{Vaz16b} focused on this, included further system classes and described two basic concepts to deal with bilateral control. The first \cross{concept} is to adjust the applied integral transformation leading to a symmetric backstepping transformation, which was applied in \cite{Vaz16c}.
This approach is also utilized in \cite{Bek19} for the bilateral stabilization and observer design of a viscous Hamilton-Jacobi PDE. By the help of a Hopf-Cole transformation, this semilinear parabolic system is mapped into a linear representation so that the results of \cite{Vaz16c} can be applied.
The other concept proposed in \cite{Vaz16b} is to reformulate the problem by suitable transformations to achieve unilateral actuation. To this end, the spatial domain is \emph{folded} at a folding point inside the domain, with the aim to get a system representation which is  similar to the classic unilateral problem so that existing results can be utilized. 
This is applied in \cite{Fre17} for the observer design of a parabolic PDE with constant coefficients, where measurements are available at both boundaries and in-domain, which allows to design two independent observers for the different parts of the domain. For hyperbolic $2\times 2$ systems with spatially varying coefficients, \cite{Wil19} considers the minimum time observer design. 

As far as the controller design of parabolic systems is concerned, the folding approach leads to new challenges. In particular, folding a scalar parabolic system leads to a $2\times 2$ system coupled via special folding boundary conditions (BCs), which need to be introduced to ensure the regularity of the solution.
In \cite{Vaz16b} the folding point is chosen to be the centre of the spatial domain. With constant parameters, this leads to equal diffusion coefficients in the folded system. This was extended by \cite{Che19,Che19b} to the case where the folding point can be chosen arbitrarily as a design parameter to freely distribute the control effort between the available boundary inputs. 
However, this leads to different diffusion coefficients in the folded system, requiring the introduction of coupling terms, the so-called well-posedness terms, in the target system to ensure the well-posedness of the kernel equations. This challenge is tackled by introducing a second integral transformation of Fredholm-type to eliminate the couplings and to obtain a target system, whose stability is easy to show.

So far, however, there exist no results for the bilateral control of coupled parabolic PDEs with spatially varying coefficients and distinct diffusion coefficients. For this setup with $n$ states, the folding approach leads to the difficult problem to stabilize a set of $2n$ coupled parabolic PDEs subject to folding BCs.

\subsection{Contribution}
In this paper, a stabilizing state feedback controller is designed for linear \emph{coupled parabolic PDEs} with bilateral actuation where all coefficients are spatially varying and the diffusion coefficients are distinct.

The design procedure is based on the folding approach, leading to a system representation with unilateral actuation. 
Hence, the usual backstepping transformation is applicable to determine the controller.
The corresponding kernel equations, however, attain an additional coupling due to the folding BCs, 
which leads to new problems for their solution. This\cross{, in particular,} becomes very challenging for coupled parabolic PDEs.

In contrast to the existing results in \cite{Che19,Che19b}, the solution of the kernel equations is based on the approach in \cite{Deu17,Ker19}. 
In particular, integral equations are directly derived for the second-order kernel equations by mapping the appearing wave equation operator into its canonical form. Then, integral equations result from formal integrations.
It is shown that despite the folding BCs, the kernel equations can be converted into proper integral equations by extending the method presented in \cite{Deu17} to the considered case of bilateral actuation. This\cross{, in particular,} requires new tricky ideas to setup the integral equations, in order to take the folding BCs into account. The resulting integral equations then have similar terms as in \cite{Deu17}, and can thus be solved by the method of successive approximations. 

Since the stability analysis of the target system resulting from the first backstepping transformation of the backstepping transformation is hindered by the coupling BCs, a Volterra-Fredholm transformation is applied to decouple the original target system into a cascade of stable parabolic systems, allowing a simple proof of stability relying on existing results. This\cross{, in particular,} allows the explicit specification of the stability margin for the closed-loop system. 
This transformation is the generalization of the approach in \cite{Che19b}. There, the Volterra part of the transformation can be significantly simplified because it involves only one state. In the general case of $n$ coupled states, the transformation leads to a novel, coupled set of boundary value problems (BVP) for the Volterra and Fredholm kernels of the decoupling transformation. With a similar approach as for the backstepping kernel, this BVP can be converted into integral equations and solved by the method of successive approximations.
Hence, the paper provides a systematic approach to deal with the bilateral backstepping control problem for a large class of coupled parabolic PDEs.

\subsection{Organisation and notation}
The considered bilateral stabilization problem is introduced in the next section. Section \ref{sec:stateFeedbackController} presents the design of the state feedback controller. Therein, the folding transformation is followed by the backstepping transformation and the decoupling transformation along with the resulting state feedback controller and the analysis of the closed-loop stability. In Section \ref{sec:solKernelEquations}, the solution of the kernel equations for the transformations into the target system is determined by converting them into proper integral equations and applying a successive approximation.
The comparison of the bilateral controller with an unilateral controller for an unstable system of two coupled parabolic PDEs confirms the effectiveness of the proposed design method in Section \ref{sec:simulation}.

\emph{Notation:}
In the paper, 
the notations $ f_z(z,z) = f_z(z,\zeta)|_{\zeta=z}$, $f_{\zeta}(z,z) = f_{\zeta}(z,\zeta)|_{\zeta=z}$ 
and $[\,\cdot\,]_{\ast}$ simplify the presentation of the results. The latter means that the expression in the square brackets has to be considered only if the condition $\ast$ is fulfilled.
Finally, variables with index $(\cdot)_i$ or $(\cdot)_{ij}$ are represented without index but boldface, \eg, $\K = K_{ij}$ if convenient.

\section{Problem formulation}
Consider the system
\begin{subequations}
	\label{eq:origSys}
	\begin{align}
		w_t(\y,t) &= \check{\Lambda}(\y) w_{\y\y}(\y,t) + \check{A}(\y) w(\y,t) \label{eq:origPDE}\\
		w_{\y}(0,t) &= B_0w(0,t) + u_0(t),	\quad t>0 \label{eq:origBC1}\\
		w_{\y}(1,t) &= B_1w(1,t) + u_1(t),	\quad t>0 \label{eq:origBC2}
	\end{align}
\end{subequations}
consisting of $n$ \emph{coupled linear parabolic PDEs} \eqref{eq:origPDE} for the state $w(\y,t) \in \mathbb{R}^{n}$, $n \geq 1$, defined on $(\y,t) \in (0,1)\times \mathbb{R}^{+}$ with the \emph{initial condition} (IC) $w(\y,0) = w_0(\y)$. The \emph{boundary conditions}  (BCs) \eqref{eq:origBC1}, \eqref{eq:origBC2} with  $B_0, B_1 \in \mathbb{R}^{n\times n}$ contain the inputs $u_0(t), u_1(t) \in \mathbb{R}^{n}$. The matrix $\check{\Lambda}(\y)$ is considered as \emph{diffusion matrix} and assumed to have a diagonal structure, \ie\ $\check{\Lambda}(\y) = \operatorname{diag}(\check{\lambda}_1(\y),\ldots,\check{\lambda}_{n}(\y))$ with the diffusion coefficients $\check{\lambda}_i \in C^{2}[0,1]$, $i=1,\ldots,n$. For a simpler presentation, they are assumed to be distinct and sorted, \ie\ $\bar{\check{\lambda}} \geq \check{\lambda}_1(\y) > \ldots > \check{\lambda}_n(\y) \geq \underline{\check{\lambda}}$, $\y \in [0,1]$, with positive constants $\bar{\check{\lambda}}$ and $\underline{\check{\lambda}}$. 
The \emph{reaction matrix} $\check{A} = [\check{A}_{ij}] \in  (C[0,1])^{n\times n}$ 
describes the coupling between the states and is the \mrk{main} source of possible instability. 

The system \eqref{eq:origSys} may result from the more general PDE
$w_t(\y,t) = (\check{\Lambda}(\y) w_{\y}(\y,t))_{\y} + \Phi(\y) w_{\y}(\y,t) + \check{A}(\y)w(\y,t)$ with a diagonal \emph{advection matrix} $\Phi \in (C^1[0,1])^{n \times n}$ after a Hopf-Cole-transformation (see \cite{Deu17}). 

In this paper, a {static} \emph{state feedback controller }
	\begin{align}
		\begin{bmatrix}
			u_0(t) \\
			u_1(t)
		\end{bmatrix}
		 = \mathcal{K}[w(t)]
		 \label{eq:stateFeedbackFormal}
	\end{align}
with the formal feedback operator $\mathcal{K}$ is determined by making use of the backstepping method. This controller has to ensure that the closed-loop system is exponentially stable with an \mrk{prescribed} rate of convergence.

\section{State feedback controller design}
\label{sec:stateFeedbackController}
\subsection{Folding transformation}
As suggested by \cite{Vaz16b,Che19,Che19b}, a folding transformation is applied to get a system representation with one-sided actuation to be able to apply backstepping with the usual Volterra-type integral transformation.
The spatial domain of the system is \emph{folded} at the \emph{folding point} $\yn \in (0,1)$, which is a design parameter. This leads to the new spatial domain $z \in [0,1]$ with
\begin{align}
	\label{eq:z(y)}
	z = f(\y) = \begin{cases}
		(\yn - \y)/\yn, & \y < \yn \\
		(\y - \yn)/(1-\yn), & \y \geq \yn
	\end{cases}
\end{align}
and the new state $x(z,t) = \operatorname{col}(\xl(z,t),\xr(z,t)) \in \mathbb{R}^{2n}$,
describing the left part $\xl(z,t) \in \mathbb{R}^n$ and the right part $\xr(z,t)\in \mathbb{R}^n$ \wrt\ $\yn$,
where
\begin{align}
	x(z,t) 
	= \begin{bmatrix}
		x_1(z,t) \\ \vdots \\ x_{2n}(z,t)
	\end{bmatrix} 
	= \begin{bmatrix}
		w_1(\yn-\yn z,t) \\
		\vdots \\
		w_n(\yn-\yn z,t) \\
		w_{1}(\yn + (1-\yn) z,t) \\
		\vdots \\
		w_{n}(\yn + (1-\yn) z,t)
	\end{bmatrix}.
	\label{eq:foldState}
\end{align}
The dynamics of the folded state reads
	\begin{subequations}\label{eq:foldedSystem}
		\begin{align}
			x_t(z,t) &= \Lambda(z) x_{zz}(z,t) + A(z) x(z,t) \label{eq:foldedPDE}\\
			\xr(0,t) &= \xl(0,t) \label{eq:foldBC1}\\
			\xl_z(0,t) &= -\underbrace{\frac{\yn}{1-\yn}}_{\symbDef{\tilde{z}_0}{ytn}} \xr_z(0,t) \label{eq:foldBC2} \\
			x_z(1,t) &= \underbrace{\begin{bmatrix}
				-\yn (B_0 w(0,t) + u_0(t)) \\
				(1-\yn) (B_1 w(1,t) + u_1(t))
			\end{bmatrix}}_{u(t)},
			\label{eq:foldBC3}
		\end{align}
\end{subequations}
	in which $\Lambda(z) = \operatorname{bdiag}(\Lambda_{l}(z),\Lambda_{r}(z)) =  \operatorname{diag}(\lambda_1(z),\ldots,\lambda_{2n}(z))$ with
	\begin{align}\label{eq:lambdaFolding}
		\lambda_i(z) = \begin{cases}
			\check{\lambda}_i(\yn - \yn z)/{\yn^{2}}, & i \leq n \\
			\check{\lambda}_i(\yn + (1-\yn)z)/{(1-\yn)^{2}}, & i > n,
		\end{cases}
	\end{align}
	$A(z)= [A_{ij}(z)] = \operatorname{bdiag}(A_l(z),A_r(z)) \in \mathbb{R}^{2n\times 2n}$, $A_l(z), A_r(z) \in \mathbb{R}^{n\times n}$ and%
	\begin{subequations}
		\begin{align}
			A_{l}(z) &= \check{A}(\yn - \yn z), \\
			A_{r}(z) &= \check{A}(\yn + (1- \yn) z)
		\end{align}
	\end{subequations}
	results from inserting $\y = f^{-1}(z)$ (see \eqref{eq:z(y)}) in $\check{\Lambda}(\y)$ and $\check{A}(\y)$.
Obviously, the folded system \eqref{eq:foldedSystem} is subject to the usual unilateral control but contains the \emph{folding BCs} \eqref{eq:foldBC1} and \eqref{eq:foldBC2}. They ensure that the solution of the folded system has the same regularity properties as the original system, which requires continuity of the state and its first spatial derivative at the folding point. The appearance of these BCs is the main reason for new challenges in the backstepping control design for \eqref{eq:origSys}.

To be able to apply the backstepping method presented in \cite{Deu17,Ker19}, the diffusion coefficients $\lambda_i(z)$ resulting from the folding transformation must differ at each point $z$. While the restriction to distinct diffusion coefficients is {considered} for the sake of a simpler presentation and can be extended to allow equal coefficients without obstacles, the case of intersecting coefficients, \ie\ being equal only on a finite subset of the spatial domain, has not yet been considered.
Consequently, depending on the shape of the diffusion coefficients, the folding point must be chosen appropriately to ensure $\bar{\lambda} \geq \lambda_1(z) > \ldots > \lambda_{2n}(z) \geq \underline{\lambda} > 0$, $z \in [0,1]$. However, this is no strong restriction, since due to the different scaling with $\frac{1}{\yn^2}$ and $\frac{1}{(1-\yn)^2}$ according to \eqref{eq:lambdaFolding}, appropriate folding points $\yn$ can always be found to ensure this property. {This result is valid for} both small enough and large enough folding points, where the latter simply requires a reordering of the folded states.
Moreover, in the included special case of constant diffusion coefficients, the folding point can be chosen without any restriction.

In the next steps, the state feedback controller \eqref{eq:stateFeedbackFormal} is calculated by designing an intermediate feedback law for the new input $u(t) \in \mathbb{R}^{2n}$ in \eqref{eq:foldBC3} using the backstepping method.

\subsection{Backstepping transformation}
Following the procedure known from the unilateral case, the boundedly invertible \emph{backstepping transformation} 
\begin{equation}
	\xt(z,t) = x(z,t) - \int_{0}^{z} K(z,\zeta) x(\zeta,t)\d\zeta,
	\label{eq:bsTrafo}
\end{equation}
into the new state $\xt(z,t) = \operatorname{col}(\xtl(z,t), \xtr(z,t))$
and the state feedback
\begin{align}\label{eq:controlLawFolded}
	u(t) = \tilde{u}(t) + K(1,1) x(1,t) + \int_0^{1}K_{z}(1,\zeta)x(\zeta,t)\d\zeta
\end{align}
with the kernel $K(z,\zeta) \in \mathbb{R}^{2n\times 2n}$ and the new input $\tilde{u}(t) = \operatorname{col}(0,\tilde{u}^{r}(t)) \in \mathbb{R}^{2n}$ with $\tilde{u}^{r}(t) \in \mathbb{R}^{n}$
are applied to map the system \eqref{eq:foldedSystem} into the \emph{intermediate target system}%
\begin{subequations}\label{eq:targetSys}
	\begin{align}
		{\xt}_t(z,t) =&\ \Lambda(z) \xt_{zz}(z,t) - \mu \xt(z,t)  \notag \\& - \At_0(z) \xt(0,t) - \At_1(z) \xt_z(0,t) \label{eq:targetPDE}\\
		\xtr(0,t) =&\ \xtl(0,t) \label{eq:targetBC1}\\
		\xtl_z(0,t) =&\ -\ytn\xtr_z(0,t) \label{eq:targetBC2}\\
		\xt_z(1,t) =&\ \tilde{u}(t), \label{eq:targetBC3}
	\end{align}
\end{subequations}
with \eqref{eq:targetPDE} defined on $(z,t) \in (0,1)\times \mathbb{R}^{+}$,
where $\mu \in \mathbb{R}$ is the design parameter to specify the stability margin of the system. The matrices $\At_0(z) = [\At_{0,ij}(z)] \in \mathbb{R}^{2n\times2n}$, $\At_1(z) = [\At_{1,ij}(z)] \in \mathbb{R}^{2n\times2n}$ are strictly lower triangular and are introduced to ensure well-posedness of the kernel equations for $K(z,\zeta)$ (cf. \cite{Vaz16a,Deu17}). Representing them by
	\begin{align}\label{eq:couplingMatrices}
		\At_i(z) = \begin{bmatrix}
			\At_i^{l}(z) & 0 \\
			\At_i^{lr}(z) & \At_i^{r}(z)
		\end{bmatrix},
	\end{align}
$i = 0,1$, shows that they lead to a coupling between the left and right subsystem via $\At_i^{lr}(z)$, which will be eliminated in the next step. To simplify the presentation in the following chapters, note that by inserting \eqref{eq:targetBC1} and \eqref{eq:targetBC2} into the coupling terms in \eqref{eq:targetPDE}, 
the PDE \eqref{eq:targetPDE}
can be rewritten in the form
$
	\xt_t(z,t) =\ \Lambda(z) \xt_{zz}(z,t) - \mu\xt(z,t)
	- \Ab_{0}(z) \xt(0,t) - \Ab_{1}(z) \xt_z(0,t)
$
with
\begin{align}\label{eq:couplingMatricesNew}
	\bar{A}_i(z) = \begin{bmatrix}
		\At_i^{l}(z) & 0 \\
		\bar{A}_i^{lr}(z) & 0
	\end{bmatrix}, \quad i=0,1.
\end{align}

\subsection{Decoupling transformation}
Due to the coupling BCs and the well-posed matrices $\At_0(z)$, $\At_1(z)$, the intermediate target system \eqref{eq:targetSys} has an involved structure, which hinders the stability analysis. To facilitate this problem, the boundedly invertible \textit{Volterra-Fredholm type decoupling transformation}
\begin{align}
	\xfr (z,t) = \xtr(z,t) &- \int_0^{z} Q(z,\zeta) \xtr(\zeta,t)\d\zeta \notag \\ 
	&
 	- \int_0^{1}P(z,\zeta)\xtl(\zeta,t)\d\zeta \label{eq:decouplingTrafo}
\end{align}
with the new state $\xfr(z,t) \in \mathbb{R}^{n}$ and the kernels $P(z,\zeta)$, $Q(z,\zeta) \in \mathbb{R}^{n\times n}$ is applied to the $\xtr$-system and the additional
state feedback
\begin{align}\label{eq:controlLawAddition}
\tilde{u}^{r}(t) =  \int_0^{1}
		\tilde{R}_f(\zeta)
		\xt(\zeta,t)\d\zeta
\end{align}
with $\tilde{R}_f(\zeta) = [
P_z(1,\zeta) \ Q_{z}(1,\zeta)] \in \mathbb{R}^{n \times 2n}$
is utilized. They map \eqref{eq:targetSys} into the \emph{final target system}
\begin{subequations}\label{eq:targetSysFinal}
	\begin{align}
		{\xt}^{l}_t (z,t) &= \Lambda_{l}(z)\xtl_{zz}(z,t) - \mu \xtl(z,t) \notag \\ &- \At_0^{l}(z) \xtl(0,t) - \At_1^{l}(z) \xtl_z(0,t) \\
		{\xf}^{r}_t(z,t) &= \Lambda_{r}(z)\xfr_{zz}(z,t) - \mu \xfr(z,t) \notag \\ &- \Af_0^{r}(z) \xfr(0,t) - \Af_1^{r}(z) \xfr_z(0,t) \\
		\xfr(0,t) &= \xtl(0,t) \\
		\xtl_z(0,t) &= -\ytn \xfr_z(0,t) \\
		\xtl_z(1,t) &= \xfr_z(1,t) = 0,
	\end{align}
\end{subequations}
with the strictly lower triangular coupling matrices $\Af_i^{r}(z) \in \mathbb{R}^{n\times n}$, $i=0,1$,
in which the coupling of the right and left states is removed in the PDEs. This will result in a cascade of parabolic PDEs after unfolding.

\mrk[2]{To show that \eqref{eq:decouplingTrafo} is \mrk[3]{boundedly} invertible, write \mrk{it} as
\begin{align}
\xfr(z,t) = \T_v[\xtr(t)](z) - \int_0^{1}P(z,\zeta) \xtl(\zeta,t)\d\zeta,
\end{align}
where $\T_v[\xtr(t)](z) = \xtr(z,t) - \int_0^{z} Q(z,\zeta) \xtr(\zeta,t)\d\zeta$ is a Volterra-type integral operator. Hence, the inverse transformation reads
\begin{align}
\label{eq:invTrafo}
\xtr(z,t) = \T_v^{-1}\Bigl[\xfr(\cdot,t) + \int_0^{1}P(\cdot,\zeta)\xtl(\zeta,t)\d\zeta\Bigr](z),\!
\end{align}
which requires solely the inversion of the Volterra-type integral operator $\T_v$. This inverse inherently exists and is bounded with the same reasoning as for \eqref{eq:bsTrafo}.} Therefore, the inverse transformation \eqref{eq:invTrafo} is bounded.

To be able to realize \eqref{eq:controlLawAddition}, the backstepping transformation \eqref{eq:bsTrafo} needs to be inserted for the target system state $\xt(z,t)$. After changing the order of integration, the additional state feedback then reads
\begin{align}
	\tilde{u}^{r}(t) = \int_0^{1}\underbrace{\Big( \tilde{R}_f(\zeta) - \int_{\zeta}^{1}\tilde{R}_f(z)K(z,\zeta)\d z\Big)}_{\check{R}_f(\zeta)} x(\zeta,t)\d \zeta
\end{align}
and after inserting into \eqref{eq:controlLawFolded}
\begin{multline}
	u(t) = K(1,1) x(1,t) \\ + \int_0^{1}\underbrace{\Big(K_{z}(1,\zeta) + \operatorname{col}\bigl(0,\check{R}_f(\zeta)\bigr) \Big)}_{R_f(\zeta)}x(\zeta,t)\d\zeta
	\label{eq:finalStateFeedbackFolded}
\end{multline}
is the final state feedback in terms of the folded state $x(z,t)$.

With \eqref{eq:finalStateFeedbackFolded} and the folding definitions \eqref{eq:foldBC3}, \eqref{eq:foldState}, the state feedback controller \eqref{eq:stateFeedbackFormal} for the original system directly follows as%
\begin{subequations}\label{eq:controlLaw}
	\begin{align}
	\begin{bmatrix}
	u_0(t) \\
	u_1(t)
	\end{bmatrix}
	=&\
	- \begin{bmatrix}
	B_0 w(0,t) \\ B_1 w(1,t)
	\end{bmatrix}  +
	\begin{bmatrix}
	-\frac{1}{\yn}I_n & 0 \\0 & \frac{1}{1-\yn}I_n
	\end{bmatrix}  \\&
	\cdot \Bigl(
	R_0 w(0,t) + R_1 w(1,t) + \int_0^1 R(\check{\zeta})w(\check{\zeta},t)\d\check{\zeta}
	\Bigr) \notag
	\end{align}
	with $R_b = [R_{b,ij}] \in \mathbb{R}^{2n \times n}$, $b = 0,1$, $R(\check{\zeta}) = [R_{ij}(\check{\zeta})]\in \mathbb{R}^{2n\times n}$ and 
	\begin{align}
	R_{0,ij} &= K_{ij}(1,1) \label{eq:feedbackGain1}\\
	R_{1,ij} &= K_{ij+n}(1,1) \\
	R_{ij}(\check{\zeta}) &= \begin{cases}
	\frac{1}{\yn}R_{f,ij}(\frac{\yn-\check{\zeta}}{\yn}), & \check{\zeta} < \yn \\
	\frac{1}{1-\yn}R_{f,ij+n}(\frac{\check{\zeta}-\yn}{1-\yn}), & \check{\zeta} \geq \yn,
	\end{cases}
	\label{eq:feedbackGain3}
	\end{align}
	$i=1,\ldots,2n$, $j=1,\ldots,n$.
\end{subequations}

\subsection{Stability of the target system}
To analyse the stability of the target system, \eqref{eq:targetSysFinal} is unfolded back to the original representation with 
\begin{align}
	\wt(\y,t) = 
	\begin{cases}
		\xtl((\yn - \y)/\yn,t), & \y < \yn \\
		\xfr((\y-\yn)/(1-\yn),t),& \y \geq \yn
	\end{cases}
\end{align}
leading to 
\begin{subequations}\label{eq:unfoldedTargetSystem}
	\begin{align}
		\wt_t(\y,t) &= \check{\Lambda}(\y) \wt_{\y\y}(\y,t) - \mu \wt(\y,t) \notag \\ 
		&
		- A_{0}(\y)\wt(\yn,t) - A_1(\y) \wt_{\y}(\yn,t) \label{eq:targetUnfoldPDE}\\
		\wt_{\y}(0,t) &= 0 \label{eq:targetUnfoldBC1}\\
		\wt_{\y}(1,t) &= 0 \label{eq:targetUnfoldBC2},
	\end{align}
\end{subequations}
where $A_b(\y) = [A_{b,ij}(\y)]$, $b=0,1$, $i,j=1,\ldots,n$, with
\begin{subequations}
	\begin{align}
		A_{0,ij}(\y) &= \begin{cases}
			\At^{l}_{0,ij}(z), & \y < \yn \\
			\Af^{r}_{0,ij}(z), & \y \geq \yn 			
		\end{cases}\!\!\\
		A_{1,ij}(\y) &= \begin{cases}
			-\yn\At^{l}_{1,ij}(z), & \y < \yn \\
			(1-\yn)\Af^{r}_{1,ij}(z), & \y \geq \yn 
		\end{cases}
	\end{align}
\end{subequations}
are strictly lower triangular matrices when taking \eqref{eq:z(y)} into account. Therefore, the unfolded target system \eqref{eq:unfoldedTargetSystem} has the well-known form, \ie\ it is a cascade of exponentially stable parabolic PDEs (see \cite{Vaz16,Deu17}), except for the fact that the coupling due to $A_i(\y)$, $i=0,1$, occurs at the folding point $\yn$ rather than at the left boundary. This, however, 
does not change the reasoning in the corresponding stability proof of \cite{Deu18} so that the following \mrk[2]{lemma} is \mrk{valid}.
\begin{lem}[Stability of the target system]%
	\label{thm:stabTarget}~\\%
	Assume that $\mu > 0$.
	Then, the initial value problem (IVP) for the unfolded target system \eqref{eq:unfoldedTargetSystem} is well-posed and \eqref{eq:unfoldedTargetSystem} is exponentially stable in the weighted $L_2$-norm $\|h\| = (\int_0^1\|\check{\Lambda}^{-\frac{1}{2}}(\y)h(\y)\|^2_{\mathbb{C}^n}\d \y)^{1/2}$, i.\:e.,
	\begin{equation}\label{expstab}%
	\|\wt(t)\| \leq \widetilde{M}\mathrm{e}^{(-\mu+c)t} \|\wt(0)\|,\quad t \geq 0
	\end{equation}%
	for all $\wt(0) \in (L_2(0,1))^n$ an $\widetilde{M} \geq 1$ and any $c>0$ such that $-\mu+c<0$.
\end{lem}
\begin{pf}
	By noting that the point evaluation is a relatively bounded operator for all $\yn\in (0,1)$ (see \cite[Ch. IV, \S1, Sec. 2]{Kat95}), the proof directly follows from \cite{Deu18}.
		Since the latter result ensures that the corresponding system operator is the infinitesimal generator of an analytic $C_0$-semigroup, choosing the ICs in $L_2(0,1)$ leads to a unique mild solution of \eqref{eq:unfoldedTargetSystem} verifying well-posedness (see \cite[Lem 3.1.5]{Cu95}).
		Furthermore, the decay rate is determined by the spectrum of the system operator (for details, see \cite{Deu18}).
	 \hfill$\square$
\end{pf}
Note that the decoupling of the left and right subsystems by \eqref{eq:decouplingTrafo} leads to a significant simplification of the target system structure so that the stability margin can be explicitly assigned, like in the unilateral case.

\subsection{Closed-loop stability}
\mrk[2]{Due to the bounded invertibility of \eqref{eq:bsTrafo} and \eqref{eq:decouplingTrafo}, 
the stability of the unfolded target system \eqref{eq:unfoldedTargetSystem}} according to \mrk[2]{Lemma} \ref{thm:stabTarget} implies the stability of the closed-loop folded system \eqref{eq:foldedSystem} with \eqref{eq:controlLawFolded}. After unfolding back to the system representation of \eqref{eq:origSys}, this directly leads to the following theorem.
\begin{thm}[Closed-loop stability]%
	\label{thm:closedLoop}%
	Assume that $\mu>0$. 
	Then, the IVP for the closed-loop system \eqref{eq:origSys} with \eqref{eq:controlLaw} is well-posed and the system is exponentially stable in the weighted $L_2$-norm $\|h\| = (\int_0^1\|\check{\Lambda}^{-\frac{1}{2}}(\y)h(\y)\|^2_{\mathbb{C}^n}\d \y)^{1/2}$, \ie\ 
		\begin{equation}\label{expstabClosedLoop}%
		\|{w}(t)\| \leq {M}\e^{(-\mu+c)t}\|{w}(0)\|,\quad t \geq 0
		\end{equation}%
		for all ${w}(0) \in (L_2(0,1))^n$, an ${M} \geq 1$ and any $c>0$ such that $-\mu+c<0$.
\end{thm}

\section{Solution of the kernel equations}
\label{sec:solKernelEquations}
To allow the calculation of the feedback gains \eqref{eq:feedbackGain1}--\eqref{eq:feedbackGain3}, the kernels $K(z,\zeta)$, $Q(z,\zeta)$ and $P(z,\zeta)$ of both the backstepping transformation \eqref{eq:bsTrafo} and the decoupling transformation \eqref{eq:decouplingTrafo} need to be determined.
\subsection{Backstepping transformation}
In order to map \eqref{eq:foldedSystem} into the intermediate target system \eqref{eq:targetSys}, $K(z,\zeta)$ must be the solution of the \emph{kernel equations}%
\begin{subequations}\label{eq:kernelEquations}%
\begin{flalign}%
	\label{eq:KPDE}
	&\Lambda(z) K_{zz}(z,\zeta) - (K(z,\zeta)\Lambda(\zeta))_{\zeta\zeta} = K(z,\zeta)(A(\zeta)+\mu I) \\
	&\Lambda(z) K'(z,z) + \Lambda(z) K_{z}(z,z) + K_{\zeta}(z,z)\Lambda(z)  \notag \\ 
	&\quad 
		+ K(z,z)\Lambda'(z) = -(A(z)+\mu I) \\
	&K(z,z)\Lambda(z) -\Lambda(z)K(z,z) = 0 \\
	& K(0,0) = 0 \\
	& K(z,0) \Lambda(0) S_1 + \At_1(z) S_1 = 0 \label{eq:kernelFoldBC1}\\\
	& K_{\zeta}(z,0)\Lambda(0) S_2 + K(z,0)\Lambda'(0) S_2 - \At_0(z) S_2 = 0 \label{eq:kernelFoldBC2} \hspace{-4mm}
\end{flalign}%
\end{subequations} %
with \eqref{eq:KPDE} defined on $0< \zeta < z < 1$ and
	\begin{align}
		S_1 = \begin{bmatrix}
			-\ytn I_n \\ I_n
		\end{bmatrix}, \quad 
		S_2 = \begin{bmatrix}
			I_n \\ I_n
		\end{bmatrix},
	\end{align}
which follow from the same calculations as in \cite{Vaz16a,Deu17}. 
The kernel equations \eqref{eq:kernelEquations} are similar to the ones found in \cite{Deu17} except for the new folding BCs \eqref{eq:kernelFoldBC1}, \eqref{eq:kernelFoldBC2}, requiring significant modifications of the solution procedure. In the remainder of this section, the following theorem will be proved.
\begin{thm}[Kernel equations of the Volterra transformation]\label{thm:kernelEquations}
	The kernel equations \eqref{eq:kernelEquations} have a piecewise {continuous} solution on the spatial domain $0\leq \zeta \leq z \leq 1$.
\end{thm}%

This result, in particular, means that \eqref{eq:bsTrafo} exists and maps $L_2$-functions into $L_2$-functions, which is sufficient for the shown closed-loop stability property. Depending on the regularity of the system parameters, it is also possible to verify a higher regularity of the kernel.
%

\subsubsection{Canonical kernel equations}%
\label{sec:canonKernelBackstepping}
The solution of the kernel equations \eqref{eq:kernelEquations} relies on a transformation into integral equations and their solution with the method of successive approximations using the approach in \cite{Ker19}. In contrast to \cite{Deu17}, it is shown in \cite{Ker19} that the step of eliminating the first-order derivatives in the canonical kernel equations is not necessary, simplifying the derivation of the integral equations. 

As first step, the kernel equations \eqref{eq:kernelEquations} are considered for each matrix element $\K(z,\zeta) = K_{ij}(z,\zeta)$, $i,j = 1,\ldots,2n$.
For this, it is convenient to first  evaluate the component form of the BC \eqref{eq:kernelFoldBC1}, resulting in%
\begin{subequations}\label{eq:compFoldBC1Both}%
\begin{multline}\label{eq:compFoldBC1Deriv0}%
	\ytn\K(z,0)\lambda_j(0) + \ytn\At_{1,ij}(z) \\*= K_{ij+n}(z,0)\lambda_{j+n}(0)  + \At_{1,ij+n}(z)
\end{multline}%
for $j\leq n$ which is equivalent to%
\begin{multline}
	\label{eq:compFoldBC1Deriv}%
	\ytn K_{ij-n}(z,0)\lambda_{j-n}(0) + \ytn\At_{1,ij-n}(z) \\*= \K(z,0) \lambda_{j}(0)  + \At_{1,ij}(z)
\end{multline}%
\end{subequations}%
with $j >n$. This shows that the respective left ($j\leq n$) and right ($j>n$) elements of the kernel $K(z,\zeta)$ and the matrix $\At_1(z)$ are coupled. Similarly, the component form of \eqref{eq:kernelFoldBC2} reads
\begin{flalign}
	&\K_{\zeta}(z,0)\lambda_j(0) + \K(z,0)\lambda_j'(0) - \At_{0,ij}(z) -\At_{0,ij+n}(z) \notag \\*
		&
		\  =
		- K_{ij+n,\zeta}(z,0)\lambda_{j+n}(0) - K_{ij+n}(z,0)\lambda_{j+n}'(0). \!\!\!&
		\label{eq:compFoldBC2Deriv}
\end{flalign}
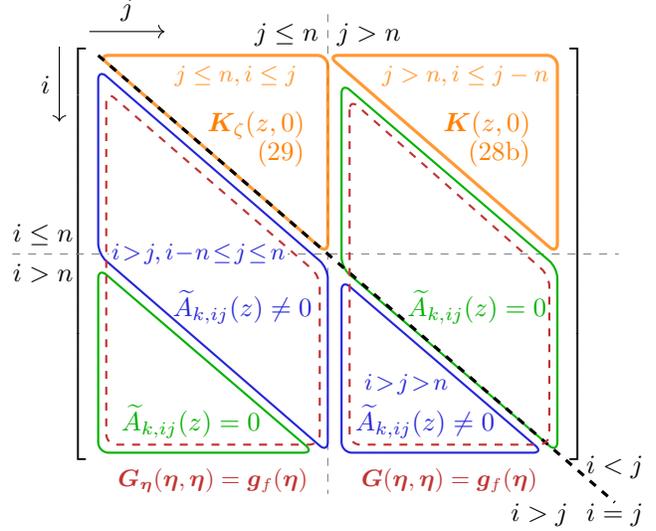
\begin{figure}%
	\colorlet{green}{greenCol}%
	\colorlet{red}{redCol}%
	\colorlet{blue}{blueCol}%
	\centering%
	\tikzset{%
add/.style args={#1 and #2}{to path={%
 ($(\tikztostart)!-#1!(\tikztotarget)$)--($(\tikztotarget)!-#2!(\tikztostart)$)%
  \tikztonodes}}%
}%
\newlength{\mywidth}
\newlength{\myheight}
\begin{tikzpicture}
	\node(m)[matrix of nodes, left delimiter = {[},right delimiter = {]},
		row sep=\linewidth/4,column sep=\linewidth/3.35, 
		every node/.append style={inner sep=0pt,minimum size=10pt}] at (0,0){
			{}&{}&{} \\ 
			{}&{}&{} \\ 
			{}&{}&{} \\ 
	};
	
	\node(i leq n)[anchor=south east, xshift=-2.1mm,inner xsep=0pt] at (m.west){$i \leq n$};
	\node[anchor=north east, xshift=-2.1mm,inner xsep=0pt] at (m.west){$i > n$};
	\node(j leq n)[anchor=south east] at (m-1-2.north){$j \leq n$};
	\node[anchor=south west] at (m-1-2.north){$j > n$};
	\getheightofnode{\myheight}{j leq n}
	\getwidthofnode{\mywidth}{i leq n}
	\draw[gray,dashed](m-1-2.north)++(0,\myheight)--
		($(m-3-2.south)+(0,-\myheight)$);

	\draw[gray,dashed](m-2-1.west-|m.west)++(-\mywidth-2.1mm,0)--($(m-2-3.east)+(\mywidth,0)$);

	\draw[rounded corners,orange,very thick,opacity = 0.8]
	(m-1-1.north west)--(m-2-2.center)--(m-1-2.north)--node[below,font=\footnotesize,orange,pos=0.4]{$j\leq n, i \leq j$}cycle;
 	\coordinate[yshift=5mm](center1) at (barycentric cs:m-1-1=1,m-2-2=1,m-1-2=1) {};
	\draw[rounded corners,orange,very thick,opacity = 0.8]
	(m-1-2.north)--(m-2-3.east)--(m-1-3.north east)--node[below,font=\footnotesize,orange,pos=0.4]{$j > n, i \leq j-n$}cycle;
	\coordinate[yshift=5mm] (center2) at (barycentric cs:m-1-2=1,m-2-3=1,m-1-3=1) {};
	
	\draw[blue,rounded corners,thick](m-1-1.west)--(m-2-1.west)--(m-3-2.south)--(m-2-2.south)--
	cycle;
	\node[font=\footnotesize,rounded corners,blue, fill=white,inner sep=1pt] at ($(m-2-1)!.4!(m-2-2)$) {$i\!>\!j, i\!-\!n\!\leq\! j \!\leq \!n$};
	\path[name path = diag](m-3-3.south)--(m-2-2.south);
	\path[name path = down](m-2-2.east)--(m-3-2.east);
	\draw[blue,rounded corners, thick,name intersections={of=diag and down}](m-3-2.south east)--(m-3-3.south)--
	(intersection-1)--cycle;
	\node[blue,font=\footnotesize,rounded corners,anchor=north east] at ($(m-3-3)!.4!(m-2-2)$) {$i\!>\!j\! >\! n$};
	
	\draw[green,rounded corners,thick](m-1-2.south east)--(m-2-2.south east)--(m-3-3.south east)--(m-2-3.south east)--cycle;
	\draw[green,rounded corners, thick,name intersections={of=diag and down}](m-3-1.south west)--(m-3-2.south west)--(m-2-1.south west)--cycle;
	
	\coordinate[xshift=-1mm,yshift=-2mm](center3) at (barycentric cs:m-2-1=1,m-3-1=1,m-3-2=1) {};
	\coordinate[yshift=-1mm](center4) at (barycentric cs:m-2-2=1,m-3-2=1,m-3-3=1) {};
	
	\node[blue,yshift=-7mm,xshift=3mm]at ($(m-2-1)!.5!(m-2-2)$){$\At_{k,ij}(z)\neq 0$};
	\node[green,yshift=-7mm,xshift=5.5mm]at ($(m-2-2)!.5!(m-2-3)$){$\At_{k,ij}(z)= 0$};
	\node[orange,yshift=-4ex,align=right] at (center1){$\K_{\zeta}(z,0)$\\ \eqref{eq:compFoldBC2Deriv}};
	\node[orange,yshift=-4ex,xshift=2mm,align=right] at (center2){$\K(z,0)$ \\ \eqref{eq:compFoldBC1Deriv}};
	\path(m-3-1.south west)--node[green,above=3mm,anchor=base,xshift=-2mm]{$\At_{k,ij}(z) = 0$}(m-3-2.south west);
	\path(m-3-2.south)--node[blue,above=3mm,anchor=base,xshift=-0.5mm] {$\At_{k,ij}(z) \neq 0$}(m-3-3.south west);
	
\draw[very thick, dashed](m-1-1.north west)--($(m-3-3.south east)!-0.125!(m-1-1.north west)$)node[anchor=north,inner sep=0pt](rightBottom){$i=j$};
\node[left= 0.1 of rightBottom]{$i>j$};
\node[above =0.2 of rightBottom]{$i<j$};

\draw[->]([yshift=2mm]m.north west)--node[above=4pt,inner sep=0pt]{$j$}++(1cm,0);
\draw[->]([xshift=-4mm]m.north west)--node[left]{$i$}++(0,-1cm);

\begin{scope}[every node/.append style={inner sep=1mm}]
	\node(m11)[fit=(m-1-1.south west)]{};
	\node(m31)[fit=(m-3-1.south west)]{};
	\node(m32)[fit=(m-3-2.south)]{};
	\node(m22)[fit=(m-2-2.south)]{};

	\node(m32east)[fit=(m-3-2.south east)]{};	
	\node(m33)[fit=(m-3-3.south east)]{};
	\node(m23)[fit=(m-2-3.south east)]{};
	\node(m12)[anchor=north]at(m-1-2.south east){};

	\draw[red,dashed,rounded corners, thick]
		(m11.south east)--(m31.north east)--node[red,below=2mm,font=\footnotesize]{$\pG_{\peta}(\peta,\peta) = \bm{g}_f(\peta)$}(m32.north west)--(m22.south west)--cycle;
	\draw[red,dashed,rounded corners, thick]
			(m12.south east)--(m32east.north east)--node[red,below=2mm,font=\footnotesize]{$\pG(\peta,\peta) = \bm{g}_f(\peta)$}(m33.north west)--(m23.south west)--cycle;
\end{scope}%
\end{tikzpicture}%
%
	\caption{Index combinations $i,j = 1,\ldots,2n$, for which the BCs \eqref{eq:compFoldBC1Both}, \eqref{eq:compFoldBC2Deriv} are either fulfilled by the kernel BCs or the choice of the coupling matrix elements $\At_{0,ij}(z)$, $\At_{1,ij}(z)$.
		The orange triangles mark the index combinations for which the kernel elements must fulfil the coupling BCs \eqref{eq:compFoldBC2Diag}, \eqref{eq:compFoldBC1} and \eqref{eq:compFoldBC2}, the blue areas represent the indices where $\At_{k,ij}(z) \neq 0$, $k=0,1$, to fulfil \eqref{eq:A01}. For the green regions, $\At_{k,ij}(z) = 0$ but \eqref{eq:compFoldBC1Both} and \eqref{eq:compFoldBC2Deriv} are fulfilled due to the appearing coupling, \ie\ by the blue counterparts \wrt\ the line $j=n$. The parts where artificial BCs \eqref{eq:wellPosedBC}, \eqref{eq:wellPosedBC2} need to be introduced are marked by the red dashed borders.}%
	\label{fig:areas}%
\end{figure}%
In the following, \eqref{eq:compFoldBC1Deriv} is considered as a BC for the right elements ($j>n$), whereas \eqref{eq:compFoldBC2Deriv} is a BC for the left elements ($j\leq n$), which is indicated by the orange areas in Figure \ref{fig:areas}.

The coupling matrices $\At_0(z)$ and $\At_1(z)$ are needed to ensure well-posedness of the kernel equations (see \cite{Deu17}), \mrk{naming them} \emph{well-posedness terms}. Their task is to fulfil the BCs \eqref{eq:compFoldBC1Both} and \eqref{eq:compFoldBC2Deriv} for some indices so that no condition on the respective kernel element results. Yet, they can only cover so much conditions that they attain a strictly lower triangular structure, \ie\ $\At_{k,ij}(z)\neq 0$, $k=0,1$, only for $i>j$ which is important for the stability of the target system. In the non-folding case, they need to remove the BC at the lower boundary $(z,0)$ for all kernel elements with $i>j$ (see \cite{Deu17}). 
In the folding case, \eqref{eq:compFoldBC1Deriv0} shows that removing the BC for elements with $j < i \leq n$ by choosing%
\begin{subequations}%
	\begin{flalign}%
		&\At_{1,ij+n}(z) = 0 \\*
		&\At_{1,ij}(z) = -\K(z,0)\lambda_j(0) + \tfrac{1}{\ytn}K_{ij+n}(z,0)\lambda_{j+n}(0)\hspace{-1cm}&
		\label{eq:At1ij}
	\end{flalign}
\end{subequations}%
for $j\leq n$
also removes the BC for the element $K_{ij+n}(z,\zeta)$.
Covering all index combinations this way and performing the same considerations for \eqref{eq:compFoldBC2Deriv} finally leads to the \emph{component form} of the kernel equations
\begin{multline}
\label{eq:compFormPDE}
\lambda_i(z)\K_{zz}(z,\zeta) - (\lambda_j(\zeta)\K(z,\zeta))_{\zeta\zeta}   \\
= \sum_{k=1}^{2n} \Big( K_{ik}(z,\zeta)A_{kj}(\zeta) \Big) + \mu\K(z,\zeta) 
\end{multline}
defined on $0 < \zeta < z < 1$, with the BCs%
\newlength{\breite}%
\settowidth{\breite}{$\K(z,z)$}%
\begin{subequations}\label{eq:compFormDiag}%
	\begin{flalign}%
	& \underline{i = j:}\ \notag &\\& \quad
	\K(z,z) = -\int_0^z\frac{A_{ii}(\zeta) + \mu}{2\sqrt{\lambda_i(\zeta)\lambda_i(z)}}\d \zeta\label{eq:KzzDiag} &\\[10pt]
	&\quad  i\leq n:\notag \\* 
	& \qquad 
	\K_{\zeta}(z,0)\lambda_i(0) + \K(z,0)\lambda_i'(0)  \label{eq:compFoldBC2Diag}&\\* 
	&
	\qquad \  = 
	- K_{ii+n,\zeta}(z,0)\lambda_{i+n}(0) - K_{ii+n}(z,0)\lambda_{i+n}'(0)\notag 
	\end{flalign}%
\end{subequations}%
and%
\begin{subequations}\label{eq:compFormNonDiag}%
	\begin{flalign}%
	&\underline{i \neq j:}\notag 
	&\\& \quad 
		\K(z,z) = 0\label{zzbedkeq}&\\
		&\quad  \K_z(z,z) = \frac{A_{ij}(z)}{\lambda_j(z) - \lambda_i(z)} \label{zzpzbed}
	&\\[10pt]  & \quad 
	j>n, \ i\leq j-n: \notag &\\*& \qquad
	\K(z,0)\lambda_j(0) = K_{ij-n}(z,0)\lambda_{j-n}(0)\ytn \label{eq:compFoldBC1}\hspace{-4mm} &\\[10pt] & \quad
	j\leq n, \ i < j: \notag &\\*&\qquad
	\K_{\zeta}(z,0)\lambda_j(0) + \K(z,0)\lambda_j'(0) \label{eq:compFoldBC2} &\\* 
	&
	\qquad \  = 
	- K_{ij+n,\zeta}(z,0)\lambda_{j+n}(0) - K_{ij+n}(z,0)\lambda_{j+n}'(0) \notag
	\end{flalign}
\end{subequations}
when the non-zero elements of $\At_0(z)$ and $\At_1(z)$ are chosen to be%
\begin{subequations}\label{eq:A01}
	\begin{flalign}
	&\underline{i>j,\ i-n \leq j \leq n:} \notag &\\*
		 & \At_{0,ij}(z)  = \lambda_j(0)\K_{\zeta}(z,0) + \lambda'_j(0) \K(z,0)
		\label{eq:A0<n}
		&\\* \notag
		& \hphantom{\At_{0,ij}(z)  =}
		+ \lambda_{j+n}(0) K_{ij+n,\zeta}(z,0)
		+ \lambda_{j+n}'(0) K_{ij+n}(z,0) &\\
		& \At_{1,ij}(z)  =
		-\lambda_j(0)\K(z,0) + \frac{1}{\ytn}\lambda_{j+n}(0)K_{ij+n}(z,0)
		\label{eq:A1<n}
		&\\
	&\underline{i>j > n:} \notag &\\* 
		&\At_{0,ij}(z)  = +\lambda_{j}(0)\K_{\zeta}(z,0) + \lambda'_j(0) \K(z,0)
		\label{eq:A0>n}
		&\\* \notag
		& \hphantom{\At_{0,ij}(z)  =}
		+ \lambda_{j-n}(0) K_{ij-n,\zeta}(z,0)
		+ \lambda_{j-n}'(0) K_{ij-n}(z,0)
		&\\
		& \At_{1,ij}(z)  =
		-\lambda_j(0)\K(z,0) + \ytn\lambda_{j-n}(0)K_{ij-n}(z,0).  &
		\label{eq:A1>n}
	\end{flalign}
\end{subequations}
Figure \ref{fig:areas} provides a graphical overview of the related index combinations for the kernel BCs and the coupling matrix elements needed to fulfil \eqref{eq:compFoldBC1Both} and \eqref{eq:compFoldBC2Deriv}.
Note that with the choice \eqref{eq:A01}, the matrices $\At_0(z)$ and $\At_1(z)$ have the structure
\begin{align}
	\At_k(z) = \begin{bmatrix}
		\At^k_{\ltriag 1}(z) & 0 \\
		\At^k_{\utriag}(z) & \At^k_{\ltriag 2}(z)
	\end{bmatrix}, \quad k=0,1,
	\label{eq:A0Struc}
\end{align}
where $\At^k_{\ltriag1}(z)$ and $\At^k_{\ltriag2}(z)$ are strictly lower triangular matrices and $\At^k_{\utriag}(z)$ is an upper triangular matrix so that in total, $\At_0(z)$ and $\At_1(z)$ are strictly lower triangular (blue regions in Figure \ref{fig:areas}).

Following the lines in \cite{Ker19}, the
component form \eqref{eq:compFormPDE} of the kernel equations is mapped into its canonical form with the new kernel elements
\begin{equation}
	\pG(\pxi,\peta) = G_{ij}(\xi_{ij},\eta_{ij}) = \lambda_j(\zeta(\pxi,\peta))K_{ij}(z(\pxi,\peta),\zeta(\pxi,\peta)).
\end{equation}
For this, the transformation to normalize the coefficients of the highest derivative to $1$ is combined with the transformation of a wave equation into its canonical form. In particular, the \textit{canonical coordinates}%
\begin{subequations}\label{eq:xieta(z,zeta)}
	\begin{align}
		\pxi &= \xi_{ij}(z,\zeta) = \tfrac{1}{2}(1-\s)(\phi_i(1)+\phi_j(1)) + 
		\s(\phi_i(z)+\phi_j(\zeta)) \\
		\peta &= \eta_{ij}(z,\zeta) = 
		-\tfrac{1}{2}(1-\s)(\phi_i(1)-\phi_j(1)) + \phi_i(z) - \phi_j(\zeta)	
	\end{align}
\end{subequations}
and their inverses%
\begin{subequations}\label{eq:zzeta(xi,eta)}
	\begin{align}
		z &= z(\pxi,\peta) =  \phi_i^{-1}(\tfrac{1}{2}(\s\pxi+\peta)+\tfrac{1}{2}(1-\s)\phi_i(1)) \\
		\zeta &= \zeta(\pxi,\peta) = \phi_j^{-1}(\tfrac{1}{2}(\s\pxi-\peta)+\tfrac{1}{2}(1-\s)\phi_j(1))
		\label{eq:zeta_inv}
	\end{align}
\end{subequations}
are considered, where
\begin{equation}\label{shdef}
	\bm{s} = s_{ij} = \begin{cases}
	1, & i \leq j\\
	-1, & i > j
	\end{cases}
\end{equation}
and 
\begin{align}
	\phi_i(z) = \int_0^z\frac{\d\oz}{\sqrt{\lambda_i(\oz)}}, \quad i=1,\ldots,2n
\end{align}
(see \cite{Ker19}).
This leads to the \emph{canonical kernel equations}%
\begin{subequations}\label{eq:canonKernelEquations}
	\begin{flalign}%
	&\bm{G}_{\bm{\xi}\bm{\eta}}(\bm{\xi},\bm{\eta})
	= \F[G,\pG_{\pxi},\pG_{\peta}](\pxi,\peta)  \notag \\
	& = \frac{\s}{4}\underbrace{\sum_{k=1}^{2n} \Big( \frac{\lambda_j(\zeta)}{\lambda_k(\zeta)}G_{ik}(\bm{\xi},\bm{\eta})A_{kj}(\zeta) \Big)\Bigr|_{\zeta=\zeta(\pxi,\peta)}}_{\bm{\mathcal{A}}[G](\pxi,\peta)} 
		+ \frac{\s\mu}{4} \pG(\bm{\xi},\bm{\eta}) &\notag\\
	& \hphantom{=\ } -\frac{1}{4} \bar{\bm{a}}_{\Delta}(\pxi,\peta) 
	\pG_{\pxi}(\bm{\xi},\bm{\eta})  
	+\frac{\s}{4}\bar{\bm{a}}_{\Sigma}(\pxi,\peta)\pG_{\peta}(\bm{\xi},\bm{\eta}) \hspace{-4mm}&\label{eq:canonPDE} 
	\end{flalign}
	with the BCs 
	\begin{flalign}
		&\underline{i=j:}\  \notag  \\& \quad
		\pG(\pxi,0) = \underbrace{-\sqrt{\lambda_i(z(\pxi,0))}\int\limits_0^{z(\pxi,0)} \frac{A_{ii}(\oz) + \mu}{2\sqrt{\lambda_i(\oz)}}\d \oz}_{\cHDef{sym:bc0}(\pxi)} \label{eq:canonBC1}\\[10pt]
		&\quad i \leq n: 
		\bigl(\pG_{\peta}(\peta,\peta)-\pG_{\pxi}(\peta,\peta)\bigr)\ytn \notag \\*
		& \quad \hphantom{i \leq n:}
		= 
		G_{ii+n,\xi}(\peta,\peta)-G_{ii+n,\eta}(\peta,\peta)
		\label{eq:canonBCFold2i=j}
		 \\[10pt]
		&\underline{i \neq j:} \notag \\ 
		& \quad \pG(\pxi,\peta_l(\pxi)) = 0 \label{eq:canonBC3}\\
		&\quad \s\pG_{\pxi}(\pxi,\peta_l(\pxi)) + \pG_{\peta}(\pxi,\peta_l(\pxi)) \notag \\
		& \quad \ =\underbrace{\frac{A_{ij}(z)\lambda_j(z)\sqrt{\lambda_i(z)}}{\lambda_j(z) - \lambda_i(z)}\Bigg|_{z=z(\pxi,\peta_l(\pxi))}}_{\cHDef{sym:bc1}(\pxi)} \label{eq:canonBC4}\\[10pt]
		&\quad j > n,\ i \leq j-n: 
			\pG(\peta,\peta)  = \ytn G_{ij-n}(\peta,\peta) \label{eq:canonBCFold1} \\[10pt]&
			\quad j \leq n,\ i < j: 
			\bigl(\pG_{\peta}(\peta,\peta)-\pG_{\pxi}(\peta,\peta)\bigr)\ytn \label{eq:canonBCFold2}\\*
			& \quad 
			\hphantom{j \leq n,\ i < j:} = 
			G_{ij+n,\xi}(\peta,\peta)-G_{ij+n,\eta}(\peta,\peta)
			\notag
		 \\[10pt]
		&\quad j > n,\ i > j-n: 
		\pG(\peta,\peta) = \bm{g}_f(\peta) \label{eq:wellPosedBC}\\[10pt]
		&\quad j \leq n,\ i > j: 
		\pG_{\peta}(\peta,\peta) = \bm{g}_f(\peta) \label{eq:wellPosedBC2}  
	\end{flalign}
\end{subequations}
and%
\begin{subequations}
	\begin{align}
		\bar{\bm{a}}_{\Delta}(\pxi,\peta) &= \tfrac{\lambda_i'(z)}{2\sqrt{\lambda_i(z)}}-
		\tfrac{\lambda_j'(\zeta)}{2\sqrt{\lambda_j(\zeta)}} \\
		\bar{\bm{a}}_{\Sigma}(\pxi,\peta) &= \tfrac{\lambda_i'(z)}{2\sqrt{\lambda_i(z)}}+
		\tfrac{\lambda_j'(\zeta)}{2\sqrt{\lambda_j(\zeta)}}
	\end{align}
\end{subequations}
where $(z,\zeta)$ are substituted by \eqref{eq:zzeta(xi,eta)}. 
In \eqref{eq:canonBCFold2i=j} and \eqref{eq:canonBCFold2},
the {result}
$
\sqrt{{\lambda_{j+n}(0)}/{\lambda_j(0)}} = \ytn = {\yn}/{(1-\yn)}
$
was {applied},
which can be derived from \eqref{eq:lambdaFolding}.
The fact that $\xi_{ij+n}(\peta,\peta)=\eta_{ij+n}(\peta,\peta)=\peta$ for $i \leq j \leq n$, respectively $\xi_{ij-n}(\peta,\peta) = \eta_{ij-n}(\peta,\peta) = \peta$ for $i \leq j-n$ is utilized in \eqref{eq:canonBCFold2i=j}, \eqref{eq:canonBCFold1} and \eqref{eq:canonBCFold2}.
In \eqref{eq:canonBC3} and \eqref{eq:canonBC4} $\peta_l(\pxi)$ is the strictly monotonically decreasing lower boundary of the canonical spatial domain (see \cite{Deu17} and Figure \ref{fig:intDirections}), which is also defined in the case $i=j$ as 
	$\peta_l(\pxi) = 0$.
The \emph{artificial BCs} \eqref{eq:wellPosedBC}, \eqref{eq:wellPosedBC2} have been introduced to fully determine the kernel and contain the degrees of freedom $\bm{g}_f = g_{ij} \in C[0,\bar{\peta}]$, $i,j=1,\ldots,2n$, where $\bar{\peta} = \phi_i(1)$ for $i\leq j$ and $\bar{\peta} = \phi_j(1)$ for $i > j$.  Their introduction as Dirichlet respectively Neumann BCs is defined by the way the equations will be converted into integral equations in the next step.
The indices for which they are provided are depicted by the red dashed areas in Figure \ref{fig:areas}.
Note that in contrast to the unilateral case, there also exist artificial BCs for $i \leq j$.

\subsubsection{Kernel integral equations}
\label{sec:kernel integral equations}
\begin{figure*}
\centering
\newcommand{\linew}{2pt} 
\newcommand{\beschbreit}{0.025}
\newcommand{\achsenfaktor}{1.4}
\newcommand{\achsenfaktorx}{1.2}


\tikzset{
	>={Stealth[length=2mm]},
	pics/coordsys/.style 2 args = {
		code={
		\newcommand{\alphaij}{#2}
		\begin{scope}[scale=1.5/\alphaij]
			\coordinate(orig) at (0,0);
			\draw[->](0,0)--(\achsenfaktorx+\achsenfaktorx*\alphaij,0)node[below right]{$\xi_{ij#1}$};
			\draw[->](0,{-\achsenfaktor+(\achsenfaktor*\alphaij)/2+(1/\alphaij)/4})--(0,\achsenfaktor*\alphaij)node[left]{$\eta_{ij#1}$};
			
			\coordinate(a) at (\alphaij,0);
			\coordinate(1pa) at (1+\alphaij,0);
			\draw[shift={(1pa)}](0,-\beschbreit)--(0,+\beschbreit)node[above]{$b^+$};	
			
			\coordinate(b) at (0,-1+\alphaij);
			\draw[shift={(b)}](\beschbreit,0)--(-\beschbreit,0)node[left]{$b^-$};	
			\coordinate(a2) at (0,\alphaij);
			\draw[shift={(a2)}](\beschbreit,0)--(-\beschbreit,0)node[left]{$\phi_i(1)$};
			\coordinate(0) at (0,0);
			\draw[shift={(0)}](\beschbreit,0)--(-\beschbreit,0)node[left]{$0$};

			\draw[line width=\linew](0,0)--node[pos=0.5,coordinate](leftMid){}(\alphaij,\alphaij)node[coordinate](top){};
			\draw[dashed](\alphaij,\alphaij)--(1+\alphaij,-1+\alphaij);
			\coordinate(newx) at (1+\alphaij,-1+\alphaij);
			
			\begin{scope}[x={(newx)},y={($(0,0)!1!90:(newx)$)}]
			\path[out=-40,in=170, line width=\linew,name path=gam2path](0,0) to
			 node[pos=0.5,coordinate](lowLeft1){}(0.5,0)node[coordinate](lowMid1){};
			\path[out=170-180,in=175, line width=\linew](0.5,0) to (1,0);
			\path[out=-40,in=170, line width=\linew, shift={(0,-0.03)}](0-0.01,0) to
			node[pos=0.65,coordinate](lowLeft2){}(0.5,0);
			\path[out=170-180,in=175, line width=\linew, shift={(0,-0.03)}](0.5,0) to 
			node[coordinate,pos=0.05](lowMid2){} (1+0.005,0);
			\path[out=-40,in=170, line width=\linew, shift={(0,-0.06)}](0-0.01,0) to
			node[pos=0.65,coordinate](lowLeft3){}(0.5,0)node[coordinate](lowMid3){}node[below,sloped]{$\peta_l(\pxi)$};
			\path[out=170-180,in=175, line width=\linew, shift={(0,-0.06)}](0.5,0) to (1+0.005,0);
			\end{scope}

		\end{scope}
		}
	},
	coordsys/.default={0.6}
}

\begin{tikzpicture}
\newcommand{\alphaA}{0.6}
\newcommand{\alphaB}{0.5}
\pic(left) at (0,0){coordsys={}{\alphaA}};
\pic(right) at (9,0){coordsys={}{\alphaB}};
\node[left=of leftorig,draw]{$j\leq n$:};
\node[above=0.25of lefttop,anchor=base west,xshift=-1cm]{$\textcolor{blueCol}{\pJ}(\peta,\peta) = \frac{1}{\ytn}                           \textcolor{blueCol}{H_{ij+n}}(\peta,\peta)$};
\node[left=of rightorig,draw]{$j > n$:};
\node[above=0.25of righttop,anchor=base west,xshift=-1cm]{$\textcolor{redCol}{\pG}(\peta,\peta) =  \ytn \textcolor{redCol}{G_{ij-n}}(\peta,\peta)$};

\begin{scope}[x={(leftnewx)},y={($(0,0)!1!90:(leftnewx)$)}]
	\draw[out=-40,in=170, line width=\linew,name path=gam2path,blueCol](0,0) to (0.5,0);
	\draw[out=170-180,in=175, line width=\linew,blueCol](0.5,0) to (1,0);
	\draw[out=-40,in=170, line width=\linew, shift={(0,-0.03)},redCol](0-0.01,0) to (0.5,0);
	\draw[out=170-180,in=175, line width=\linew, shift={(0,-0.03)},redCol](0.5,0) to (1+0.005,0);
	\draw[out=-40,in=170, line width=\linew, shift={(0,-0.06)},greenCol](0-0.01,0) to (0.5,0);
	\draw[out=170-180,in=175, line width=\linew, shift={(0,-0.06)},greenCol](0.5,0) to (1+0.005,0);	
\end{scope}
\begin{scope}[shift=(leftorig),scale=2.75]
	\draw[line width=\linew,blueCol](leftorig)--(lefttop);
	
	\draw[->,white,line width=3pt,line cap=butt,shorten <=\linew](leftlowLeft2)--++(0,0.4);
	\draw[->,redCol,thick](leftlowLeft2)--node[pos = 0.6,left,redCol,fill=white,inner sep=1.5pt,outer sep=1.5pt]{$\pG$}++(0,0.4);
	
	\draw[->,white,line width=3pt,line cap=butt,shorten <=\linew](leftlowMid3)--++(0,0.35);
	\draw[->,greenCol,thick](leftlowMid3)--node[pos = 0.55,left,greenCol,fill=white,inner sep=1.5pt,outer sep=1.5pt]{$\pH$}++(0,0.45);
	
	\draw[->,blueCol,thick](leftleftMid)--node[pos = 0.5,above,blueCol]{$\pJ$}++(0.5,0);
\end{scope}

\begin{scope}[shift=(rightorig),x={(rightnewx)},y={($(0,0)!1!90:(rightnewx)$)}]
	\draw[out=-40,in=170, line width=\linew,name path=gam2path,redCol](0,0) to (0.5,0);
	\draw[out=170-180,in=175, line width=\linew,redCol](0.5,0) to (1,0);
	\draw[out=-40,in=170, line width=\linew, shift={(0,-0.03)},blueCol](0-0.01,0) to (0.5,0);
	\draw[out=170-180,in=175, line width=\linew, shift={(0,-0.03)},blueCol](0.5,0) to (1+0.005,0);
\end{scope}
\begin{scope}[shift=(rightorig),scale=2.75]
	\draw[line width=\linew,redCol](rightorig)--(righttop);
	
	\draw[->,white,line width=3pt,line cap=butt,shorten <=\linew](rightlowLeft2)--++(0,0.4);
	\draw[->,blueCol,thick](rightlowLeft2)--node[pos = 0.45,left,blueCol]{$\pH$}++(0,0.45);
	\draw[->,redCol,thick](rightleftMid)--node[pos = 0.5,above,redCol]{$\pG$}++(0.5,0);
\end{scope}


\end{tikzpicture}
\caption{The spatial domains of the kernel equations \eqref{eq:canonKernelEquations} in the canonical coordinate systems $(\pxi,\peta)$ for $i\leq j$ with the different directions of integration in the cases $j\leq n$ (left) leading to \eqref{eq:intEqsNewLeft} and $j > n$ (right) resulting in \eqref{eq:intEqsNewRight} to account for the coupling BCs at $(\peta,\peta)$. 
The spatial domain is characterized by $b^- = \phi_i(1) - \phi_j(1)$, $b^+ = \phi_i(1)+\phi_j(1)$ and $\peta_l(\pxi)$, which is zero for $i=j$. The thick coloured lines represent the BCs and the coloured arrows the respective directions of integration for the corresponding variables. Herein, the blue and red arrows symbolize the first and second integration, respectively, which highlights the reversal of the integration order for $j\leq n$ and $j>n$. For the left elements, the BC for $\pJ(\peta,\peta)$ depends {on $H_{ij+n}(\peta,\peta)$} for which the respective integral equation can be substituted. On the contrary, for $j>n$, $\pG(\peta,\peta)$ depends on $G_{ij-n}(\peta,\peta)$ for which again the corresponding integral equation can be substituted. In this way, the kernel equations can be converted into proper integral equations.}
\label{fig:intDirections}
\end{figure*}
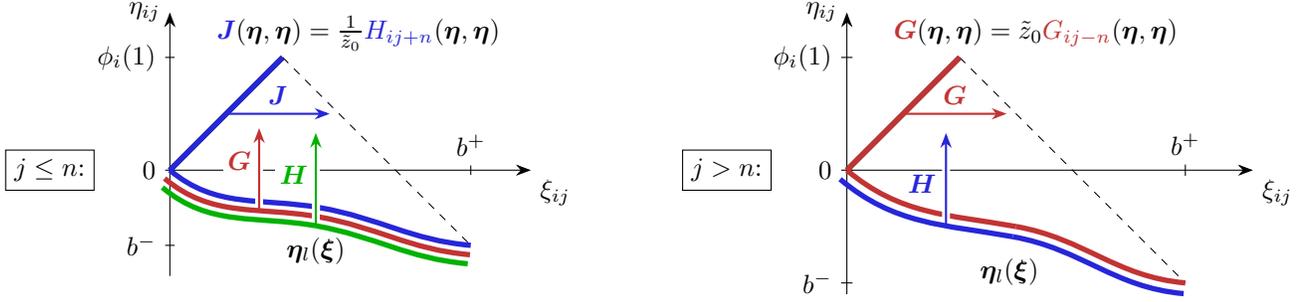
For the conversion into integral equations, 
\eqref{eq:canonPDE} is formally integrated \wrt\ $\pxi$ and $\peta$. 
Due to the different types of coupling BCs \eqref{eq:canonBCFold1} for $j>n$ and \eqref{eq:canonBCFold2i=j}, \eqref{eq:canonBCFold2} for $j\leq n$, a different sequence of the formal integration is needed for elements with $j\leq n$ and $j >n$. This is visualized in Figure \ref{fig:intDirections}.
Starting with $j > n$, \eqref{eq:canonPDE} is 
integrated \wrt\ $\peta$ first and then \wrt\ $\pxi$ as in \cite{Ker19}, leading to%
\begin{subequations}\label{eq:firstMethod}
	\begin{align}
		\label{eq:firstMethodH}
		\pG_{\pxi}(\pxi,\peta) &= \pG_{\pxi}(\pxi,\peta_l(\pxi)) + \int_{\peta_l(\pxi)}^{\peta}	 \pG_{\pxi\peta}(\pxi,\oeta) \d\oeta\\
		\pG(\pxi,\peta) &= \pG(\pxi_l(\peta),\peta) + \int_{\pxi_l(\peta)}^{\pxi} \pG_{\pxi}(\oxi,\peta)\d\oxi,
		\label{eq:firstMethodG}
	\end{align}
\end{subequations}
where%
\begin{align}\label{eq:xiL}
	\pxi_l(\peta) = \begin{cases}
		\peta, & \peta \geq 0 \\
		\peta_l^{-1}(\peta), & \peta < 0
	\end{cases}
\end{align} 
is the left boundary of the spatial domain shown in Figure \ref{fig:intDirections}. Introducing the new variable 
\begin{align}
	\pH(\pxi,\peta) &\coloneqq \pG_{\pxi}(\pxi,\peta)
\end{align}
results in%
\begin{subequations}%
	\label{eq:intEqsNewRight}%
	\begin{align}%
		\pH(\pxi,\peta) &= \pH(\pxi,\peta_l(\pxi)) + \int_{\peta_l(\pxi)}^{\peta}\F[G,\pH,\pG_{\peta}](\pxi,\oeta)\d\oeta
		\label{eq:inteqH} \\
		\label{eq:inteqG1}
		\pG(\pxi,\peta) &= \pG(\pxi_l(\peta),\peta) + \int_{\pxi_l(\peta)}^{\pxi} \pH(\oxi,\peta)\d\oxi
	\end{align}%
\end{subequations}%
for $j > n$ after $\pG_{\pxi\peta}$ has been substituted by the right side of the PDE \eqref{eq:canonPDE}.
It can be seen that \eqref{eq:intEqsNewRight}
requires suitable BCs for $\pH(\pxi,\peta_l(\pxi)) = \pG_{\pxi}(\pxi,\peta_l(\pxi))$ and $\pG(\pxi_l(\peta),\peta)$.

For $j \leq n$, \eqref{eq:canonPDE} is 
integrated in the reverse order to get%
\begin{subequations}\label{eq:secondMethod}%
	\begin{align}%
		\label{eq:secondMethodJ}%
		\pG_{\peta}(\pxi,\peta) &= \pG_{\peta}(\pxi_l(\peta),\peta) + \int_{\pxi_l(\peta)}^{\pxi} \pG_{\peta\pxi}(\oxi,\peta) \d\oxi \\
		\pG(\pxi,\peta) &= \pG(\pxi,\peta_l(\pxi)) + \int_{\peta_l(\pxi)}^{\peta} \pG_{\peta}(\pxi,\oeta)\d\oeta,
		\label{eq:secondMethodG}
	\end{align}%
\end{subequations}%
where the introduction of 
\begin{align}
	\pJ(\pxi,\peta) &\coloneqq \pG_{\peta}(\pxi,\peta)
\end{align}
provides%
\begin{subequations}
	\label{eq:intEqsNewLeft}
	\begin{align}
		& \label{eq:inteqJ} 
		\pJ(\pxi,\peta) = \pJ(\pxi_l(\peta),\peta) + \int_{\pxi_l(\peta)}^{\pxi} \F[G,\pG_{\pxi},\pJ](\oxi,\peta)\d\oxi \\
		  \label{eq:inteqG2}
		&\pG(\pxi,\peta) = \pG(\pxi,\peta_l(\pxi)) + \int_{\peta_l(\pxi)}^{\peta} \pJ(\pxi,\oeta)\d\oeta,
	\end{align}
\end{subequations}
which need BCs for $\pJ(\pxi_l(\peta),\peta) = \pG_{\peta}(\pxi_l(\peta),\peta)$ and $\pG(\pxi,\peta_l(\pxi))$.
In the following, the appearing boundary terms in \eqref{eq:intEqsNewRight} and \eqref{eq:intEqsNewLeft} are substituted by suitable BCs.

As first step, with $\peta_l(\pxi) = 0$, $i=j$, the BCs \eqref{eq:canonBC1} and \eqref{eq:canonBC3} can directly be inserted for $\pG(\pxi,\peta_l(\pxi))$ in \eqref{eq:inteqG2}.
To substitute $\pH(\pxi,\peta_l(\pxi)) = \pG_{\pxi}(\pxi,\peta_l(\pxi))$ in \eqref{eq:inteqH},
differentiate \eqref{eq:canonBC3} \wrt\ $\pxi$ to obtain
\begin{align}\label{eq:Hxietal+Jxietal}
	\pG_{\pxi}(\pxi,\peta_l(\pxi)) + \pG_{\peta}(\pxi,\peta_l(\pxi))\peta_l'(\pxi) = 0, \quad i\neq j.
\end{align}
Now, solving \eqref{eq:canonBC4} for $\pG_{\peta}(\pxi,\peta_l(\pxi))$ and inserting the result into \eqref{eq:Hxietal+Jxietal} yields
\begin{align}
	\pH(\pxi,\peta_l(\pxi)) &= \underbrace{\frac{\cH{sym:bc1}(\pxi)\peta_{l}'(\pxi)}{\s \peta_l'(\pxi)-1}}_{\cHDef{sym:bc4}(\pxi)}, \quad i\neq j, \label{eq:HxietaL}
\end{align}
where the denominator cannot be zero (see \cite{Deu17}).
Moreover, 
differentiating \eqref{eq:canonBC1} \wrt\ $\pxi$ gives
\begin{align}\label{eq:Hxieta0}
	\pH(\pxi,\peta_l(\pxi)) = \pH(\pxi,0) &= \underbrace{\d_{\pxi}\cH{sym:bc0}(\pxi)}_{\cHDef{sym:cDiff}(\pxi)}, \quad i=j.
\end{align}
The remaining boundary terms incorporate the coupling BCs \eqref{eq:canonBCFold2i=j}, \eqref{eq:canonBCFold1} and \eqref{eq:canonBCFold2} and thus need special attention.

For the formulation of $\pG(\pxi_l(\peta),\peta)$ in \eqref{eq:inteqG1}, first note that the lower boundary $\peta_l(\pxi)$ of the domain $(\pxi,\peta)$ is the left boundary of the domain for $\peta < 0$, too (see Figure \ref{fig:intDirections}), so that inserting $\pxi = \pxi_l(\peta) = \peta_l^{-1}(\peta)$ in \eqref{eq:canonBC3} yields
\begin{align}
	\pG(\pxi_l(\peta),\peta) = 0, \quad \peta < 0.
\end{align}
With \eqref{eq:xiL}, this leads to
\begin{align}
	\pG(\pxi_l(\peta),\peta) = \begin{cases}
		\pG(\peta,\peta), & \peta \geq 0 \\
		0, & \peta <0,
	\end{cases}
\end{align}
where the BCs \eqref{eq:canonBCFold1} and \eqref{eq:wellPosedBC} can be inserted for $\peta \geq 0$. Then, \eqref{eq:inteqG1} reads
\begin{align}
	\pG(\pxi,\peta) =&\ [\ytn G_{ij-n}(\peta,\peta)]_{\rightcond{i &\leq j-n \\ \peta &\geq 0}} 
	+ [\bm{g}_f(\peta)]_{\rightcond{i &> j-n\\ \peta &\geq 0}} \notag \\&
	+ \int_{\pxi_l(\peta)}^{\pxi} \pH(\oxi,\peta)\d\oxi
\end{align}
for $j>n$ (see right picture of Figure \ref{fig:intDirections}).

For $G_{ij-n}(\peta,\peta)$, the respective result \eqref{eq:inteqG2} with \eqref{eq:canonBC1} and \eqref{eq:canonBC3} is now inserted to obtain
\begin{flalign}
	& \pG(\pxi,\peta) =  
	[\bm{g}_f(\peta)]_{\rightcond{i &> j\!-\!n\\ \peta &\geq 0}}
			+ \int_{\pxi_l(\peta)}^{\pxi} \pH(\oxi,\peta)\d\oxi  \\&
		\quad + [\ytn\cH{sym:bc0}(\peta)]_{\rightcond{i=j\!-\!n}}
		+ \Bigl[\ytn\intl_{\mathclap{\eta_{l,ij-n}(\peta)}}^{\peta} J_{ij-n}(\peta,\oeta)\d\oeta
		\Bigr]_{\rightcond{i & \leq j\!-\!n \\ \peta &\geq 0}} 
		.\notag
\end{flalign}

To substitute $\pJ(\pxi_l(\peta),\peta)$ in \eqref{eq:inteqJ}, the BCs \eqref{eq:canonBC4}, \eqref{eq:canonBCFold2i=j} and \eqref{eq:canonBCFold2} need to be formulated as conditions for $\pJ = \pG_{\peta}$.
Solving \eqref{eq:canonBCFold2} for $\pG_{\pxi}(\pxi,\peta_l(\pxi))$ and inserting the result into \eqref{eq:Hxietal+Jxietal} yields
\begin{align}
	\pJ(\pxi_l(\peta),\peta) &=
		 \frac{\cH{sym:bc1}(\pxi_l(\peta))}{1-\s\peta_l'(\pxi_l(\peta))} \eqqcolon \cHDef{sym:bc5}(\peta), \quad i\neq j
		 \label{eq:JxietaL}
\end{align}
in the case $\peta < 0$ after the substitution $\pxi = \pxi_l(\peta)$. Since the denominator in \eqref{eq:HxietaL} cannot be zero, the same holds for the one in \eqref{eq:JxietaL}.

For $\peta \geq 0$, \eqref{eq:canonBCFold1} is differentiated \wrt\ $\peta$ yielding
\begin{multline}\label{eq:detaGetaeta}
	\pG_{\pxi}(\peta,\peta) + \pG_{\peta}(\peta,\peta) \\= \ytn G_{ij-n,\xi}(\peta,\peta) + \ytn G_{ij-n,\eta}(\peta,\peta)
\end{multline}
for $j>n$, $i\leq j-n$.
Now note that for an element with column index $j>n$ that is subject to \eqref{eq:detaGetaeta}, the coupled element with column $j-n$ must fulfil one of the BCs \eqref{eq:canonBCFold2i=j}, \eqref{eq:canonBCFold2}. This can be exploited by considering an index-shift $j-n \to j$ for \eqref{eq:detaGetaeta}, leading to
\begin{multline}\label{eq:detaGetaetaShifted}
	G_{ij+n\xi}(\peta,\peta) + G_{ij+n,\eta}(\peta,\peta) \\= \ytn \pG_{\pxi}(\peta,\peta) + \ytn \pG_{\peta}(\peta,\peta)
\end{multline}
for $j \leq n$, $i \leq j$.
Solving \eqref{eq:detaGetaetaShifted} for $G_{ij+n,\eta}(\peta,\peta)$ and inserting the result into \eqref{eq:canonBCFold2i=j} and \eqref{eq:canonBCFold2}
leads to 
\begin{align}
	\pJ(\peta,\peta) = 
		\tfrac{1}{\ytn} G_{ij+n,\xi}(\peta,\peta)
	\label{eq:JBCetaeta}
\end{align}%
for $i\leq j \leq n$, $\peta \geq 0$ after some rearrangements, where $G_{ij+n,\xi}(\peta,\peta) = H_{ij+n}(\peta,\peta)$ can be replaced by \eqref{eq:inteqH}.

Together with \eqref{eq:wellPosedBC2} and \eqref{eq:JxietaL}, the result \eqref{eq:JBCetaeta} can now be inserted for $\pJ(\pxi_l(\peta),\peta)$ in \eqref{eq:inteqJ}.

In \eqref{eq:inteqH} and \eqref{eq:inteqJ} both $\pG_{\pxi} = \pH$ and $\pG_{\peta} = \pJ$ appear under the integrals.
Though \eqref{eq:inteqH} is only required for $j>n$ to determine $\pG$ according to \eqref{eq:inteqG1}, it can also be utilized for $j \leq n$, because \eqref{eq:HxietaL} and \eqref{eq:Hxieta0} are valid for all $j = 1,\ldots,2n$. Therefore, $\pG_{\pxi}$ can be replaced by $\pH$ in \eqref{eq:inteqJ}.
On the contrary, the equation \eqref{eq:inteqJ} may only be used for $j\leq n$, since the utilized folding BC \eqref{eq:JBCetaeta} is only valid there. Thus, integration by parts needs to be applied in \eqref{eq:inteqH} for $j>n$ to eliminate $\pG_{\peta}$, like shown in \cite{Ker19}.

Then, substituting all considered boundary terms finally leads to 
the \emph{kernel integral equations}%
\begin{subequations}
\label{eq:intEqs}
	\begin{align}
		\pG &= \pG_0 + \F_{G}[\pH,J] \\
		\pH &= \pH_0 + \F_H[G,\pH,J] \\
		\makebox[0pt][r]{$j \leq n:$} \quad \pJ &= \pJ_0 + \F_J[G,H,J]
		\label{eq:inteqFormalJ}
	\end{align}
\end{subequations}%
with%
\begin{subequations}
\label{eq:startVals}
\begin{flalign}
	\pG_0 =&\ \left[\cH{sym:bc0}(\pxi)\right]_{\rightcond{i&= j\\ j&\leq n}}
		\!+ \left[\bm{g}_f(\peta)\right]_{\rightcond{j &> n \\ i &> j-n \\ \peta& \geq 0}}
		\!+ \left[\ytn \cH{sym:bc0}(\peta)\right]_{\rightcond{i&=j\!-\!n}} &\\[10pt]
	\pH_0 =&\ H_{0,ij} = [\tfrac{1}{4}\as(\pxi,\peta) \cH{sym:bc0}(\pxi)]_{\rightcond{j&\leq n \\ i&=j}}
		+[\cH{sym:cDiff}(\pxi)]_{\rightcond{i=j}}
		\notag &\\*&
		-[\tfrac{1}{4}\as (\pxi,0) \cH{sym:bc0}(\pxi)]_{\rightcond{i=j}} 
		+[\tfrac{1}{4}\as(\pxi,\peta) \bm{g}_f(\peta)]_{\rightcond{j&>n \\ i&>j-n}}
		\notag &\\*&
		+ [\cH{sym:bc4}(\pxi)]_{\rightcond{i&\neq j}} 
		+[\tfrac{1}{4}\as(\pxi,\peta) \ytn \cH{sym:bc0}(\peta)]_{\rightcond{i&=j-n}}
		\hspace{-4mm}&\\[10pt]
	\pJ_0 =&\ 
		[\cH{sym:bc5}(\peta)]_{\rightcond{i &\neq j \\ \peta &< 0}} 
		+ [\frac{1}{\ytn} H_{0,ij+n}(\peta,\peta)]_{\rightcond{i &\leq j \\ \peta &\geq 0}}
		+ [\bm{g}_f(\peta)]_{\rightcond{i&>j\\\peta &\geq 0}} &
\end{flalign}
\end{subequations}%
as well as%
\begin{subequations}%
\label{eq:intOps}%
\begin{flalign}%
	&\F_G[\pH,J] =
		\Bigl[ \!\intl_{\peta_l(\pxi)}^{\peta}\!\!\! \pJ(\pxi,\oeta)\d\oeta\Bigr]_{\rightcond{j&\leq n}}
		\! +\Bigl[\!\intl_{\pxi_l(\peta)}^{\pxi}\!\!\! \pH(\oxi,\peta)\d\oxi\Bigr]_{\rightcond{j&>n}}
		\notag &\\*&
		\hphantom{\F_G[\pH,J] =\ }
		+\Bigl[\ytn \!\!\!\! \intl_{{\eta_{l,ij-n}(\peta)}}^{\peta}\!\!\!\!\!\! J_{ij-n}(\peta,\oeta)\d\oeta\Bigr]
		_{\rightcond{i &\leq j-n \\ \peta&\geq 0}} \hspace{-4mm}&\\[5pt]
	&\F_H[G,\pH,J] = F_{H,ij}[G,\pH,J] = 
		\!\!\intl_{\peta_l(\pxi)}^{\peta} \!\! \Bigl(
			-\tfrac{\ad(\pxi,\oeta)}{4} \pH(\pxi,\oeta) 
			\notag &\\&
		\quad 
		+ \tfrac{\s}{4} \bigl(\tilde{\bm{a}}(\pxi,\oeta)+\mu\bigr) \pG(\pxi,\oeta) 
		+ [\tfrac{\s\as(\pxi,\peta)}{4} \pJ(\pxi,\oeta)]_{\rightcond{j&\leq n}}
		\notag &\\&
		\quad
		+ \tfrac{\s}{4}\bm{\mA}[G](\pxi,\oeta) 
				\Bigr)\d\oeta
		+\Bigl[
					\tfrac{\s \as(\pxi,\peta)}{4}\int_{\pxi_l(\peta)}^{\pxi}\pH(\oxi,\peta) \d\oxi
				\Bigr]_{\rightcond{j&>n}}
		\notag &\\&
		\quad
		+\Bigl[
					\tfrac{\as(\pxi,\peta)}{4}\ytn\int_{\eta_{l,ij-n}(\peta)}^{\peta} J_{ij-n}(\peta,\oeta)\d\oeta
				\Bigr]_{\rightcond{i &\leq j-n \\ \peta &\geq 0}} \hspace{-4mm}&
		\\[0pt]
	&\F_J[G,H,J] =
		\int_{\pxi_l(\peta)}^{\pxi}
		\Bigl(
			\tfrac{\s \as(\oxi,\peta)}{4}\pJ(\oxi,\peta) - \tfrac{\ad(\oxi,\peta)}{4}\pH(\oxi,\peta) 
			\notag &\\&
			+\tfrac{\s}{4}\bm{\mathcal{A}}[G](\oxi,\peta) + \mu \pG(\oxi,\peta)
		\Bigr)\d\oxi
		\notag &\\& 
			 -\Bigl[ \frac{1}{\ytn} F_{H,ij+n}[G,H_{ij+n},J](\peta,\peta)
		\Bigr]_{\rightcond{i &\leq j \\ \peta &\geq 0}}, \hspace{-4mm}&
\end{flalign}%
\end{subequations}%
where $\tilde{\bm{a}}(\pxi,\peta) = \partial_{\peta} \as(\pxi,\peta)$.
Note that $\F_G$ and $\F_H$ only contain $J_{ij}$ with $j\leq n$, which is determined by \eqref{eq:inteqFormalJ}.

\subsection{Decoupling transformation}
\label{sec:decouplingKernel}
After differentiating \eqref{eq:decouplingTrafo} \wrt\ time as well as using \eqref{eq:targetSys} and \eqref{eq:decouplingTrafo}, similar calculations as in \cite{Vaz16a,Deu17} show that \eqref{eq:targetSys} is mapped
into the final target system \eqref{eq:targetSysFinal} if the kernels $P(z,\zeta)$ and $Q(z,\zeta)$ are the solution of the \textit{kernel equations}%
\begin{subequations}\label{eq:decouplingBVP}
	{\setlength{\belowdisplayskip}{0pt}
	\begin{align}
	&\Lambda_{r}(z)P_{zz}(z,\zeta) - (P(z,\zeta)\Lambda_{l}(\zeta))_{\zeta\zeta} = 0 \label{eq:PPde} \\
	&P_{\zeta}(z,1)\Lambda_{l}(1) + P(z,1)\Lambda_l'(1) = 0 \\
	&P(0,\zeta) = 0 \\
	&P_{z}(0,\zeta) = 0
	\end{align}
	}{
	\setlength{\abovedisplayskip}{-.5\baselineskip}
	\begin{align}
	&\Lambda_{r}(z)Q_{zz}(z,\zeta) - (Q(z,\zeta)\Lambda_{r}(\zeta))_{\zeta\zeta} = 0 \label{eq:QPde}\\
	& Q_{\zeta}(z,z)\Lambda_{r}(z)+Q(z,z)\Lambda_{r}'(z) + \Lambda_{r}(z)Q'(z,z) \notag \\ 
	&
	\quad + \Lambda_{r}(z)Q_z(z,z) = 0 \\
	& \Lambda_{r}(z)Q(z,z) - Q(z,z)\Lambda_{r}(z) = 0 \\
	&Q(0,0) = 0
	\end{align}%
	}{%
	with the coupling BCs
	\begin{align}
		& P(z,0) \Lambda_l(0) - \tfrac{1}{\ytn} Q(z,0)\Lambda_r(0) \notag \\ 
		&
		= \Ab_1^{lr}(z) - \int_0^{z}Q(z,\zeta) \Ab_1^{lr}(\zeta) \d\zeta - \int_0^{1}P(z,\zeta) \At_1^{l}(\zeta)\d\zeta \notag \\ 
		&\quad 
		+ \tfrac{1}{\ytn}\check{A}_1^{r}(z), \label{eq:couplingBC1} \\[10pt]
		&Q_{\zeta}(z,0)\Lambda_{r}(0) + Q(z,0)\Lambda_{r}'(0) +P_{\zeta}(z,0)\Lambda_{l}(0) + P(z,0)\Lambda_{l}'(0) \notag \\ 
		&
		= -\Ab_0^{lr}(z) + \int_0^{z}Q(z,\zeta)\Ab_0^{lr}(\zeta)\d\zeta + \int_0^{1}P(z,\zeta)\At_0^{l}(\zeta)\d\zeta
		\notag \\ 
		&
		\quad + \check{A}_0^{r}(z),\label{eq:couplingBC2}
	\end{align}}%
\end{subequations}%
in which $P(z,\zeta)$ \mrk{is} defined on the rectangular domain $0\leq z \leq 1$, $0 \leq \zeta \leq 1$ and $Q(z,\zeta)$ \mrk{is} defined on the triangular domain $0 \leq \zeta \leq z \leq 1$. Since the kernel equations for $P(z,\zeta)$ and $Q(z,\zeta)$ are coupled via the BCs \eqref{eq:couplingBC1} and \eqref{eq:couplingBC2}, they cannot be solved independently, but a solution method needs to be derived allowing a simultaneous determination of \mrk{both kernels}. This challenging problem of solving coupled Fredholm-Volterra kernel equations is new in the backstepping framework and provides a general extension of the corresponding transformation presented in \cite{Che19b}.
In the sequel, the following well-posedness result for these kernel equations will be shown.
\begin{thm}[Kernel equations of the Volterra-Fredholm transformation]
	\label{thm:kernelEquationsDecoupling}
	The kernel equations \eqref{eq:decouplingBVP} have a piecewise {continuous} solution $P(z,\zeta)$ on the spatial domain $0\leq z \leq 1$, $0\leq \zeta \leq 1$ and $Q(z,\zeta)$ on the spatial domain $0\leq \zeta \leq z \leq 1$.
\end{thm}

Since both PDEs \eqref{eq:PPde} and \eqref{eq:QPde} have the same spatial differential operator as \eqref{eq:KPDE}, the approach to solve \eqref{eq:decouplingBVP} is the same as for the solution of \eqref{eq:kernelEquations}. The main difference is that \eqref{eq:PPde} is now defined on a rectangular domain and that the coupling between the kernel elements appears solely in the coupling BCs. Therefore, the component forms of the equations for both kernels \cross{$Q$ and $P$} are transformed into canonical coordinates in the next step, which can then be converted into integral equations. 
\subsubsection{Canonical kernel equations}
The equivalent \textit{component form} of \eqref{eq:decouplingBVP} for the matrix elements 
$\P(z,\zeta) = P_{ij}(z,\zeta)$, $i,j=1,\ldots,n$, reads
\begin{subequations}\label{eq:CompFormdecouplingBVPP}
	\begin{align}
	&\lambda^{r}_i(z)\P_{zz}(z,\zeta) - (\lambda_j^{l}(\zeta)\P(z,\zeta))_{\zeta\zeta} = 0 \label{eq:PCompPde} \\
	&\lambda_j^{l}(1)\P_{\zeta}(z,1) + \lambda_j^{l\prime}(1)\P(z,1) = 0 \\
	& \P(0,\zeta) = 0 \\
	& \P_{z}(0,\zeta) = 0 \label{eq:P3}
	\end{align}
\end{subequations}
and for $\Q(z,\zeta) = Q_{ij}(z,\zeta)$, $i,j=1,\ldots,n$,
\begin{subequations}\label{eq:CompFormdecouplingBVPQ}
		\begin{align}
		& \lambda_i^{r}(z)\Q_{zz}(z,\zeta) - (\lambda_j^{r}(\zeta)\Q(z,\zeta))_{\zeta\zeta} = 0 \label{eq:QCompPde}\\
		& \Q_{\zeta}(z,z)\lambda^{r}_j(z)+\Q(z,z)\lambda^{r\prime}_j(z) + \lambda_i^{r}(z)\Q'(z,z) \notag \\ 
		&
		\quad + \lambda_i^{r}(z)\Q_z(z,z) = 0 \label{eq:Q3} \\
		& (\lambda^{r}_i(z) - \lambda^{r}_j(z)) \Q(z,z) = 0 \label{eq:Q2}\\
		&\Q(0,0) = 0. \label{eq:Q4}
		\end{align}%
\end{subequations}
The coupling BCs for \eqref{eq:CompFormdecouplingBVPP} and \eqref{eq:CompFormdecouplingBVPQ} are%
\begin{subequations}\label{eq:compFormCouplingBCs}
	\begin{align}
		& \underline{i\leq j:} \notag \\ 
		&\ \lambda_j^{l}(0)\P(z,0) - \tfrac{1}{\ytn} \lambda_j^{r}(0)\Q(z,0) \notag \\ 
		&\ \quad	=  - \sum_{k=1}^{n}\int_0^{z}Q_{ik}(z,\zeta) \Ab_{1,kj}^{lr}(\zeta) \d\zeta 
		\notag \\ 
		&\	\qquad - \sum_{k=1}^{n}\int_0^{1}P_{ik}(z,\zeta) \At_{1,kj}^{l}(\zeta)\d\zeta +\Ab_{1,ij}^{lr}(z) \label{eq:couplingCompBC1}
		\\[10pt] 
		&\ \lambda_j^{r}(0)\Q_{\zeta}(z,0) + \lambda_j^{r\prime}(0)\Q(z,0) + \lambda_j^{l}(0)\P_{\zeta}(z,0) \notag \\ 
		& + \lambda_j^{l\prime}(0) \P(z,0) 
		= \sum_{k=1}^{n}\int_0^{z}Q_{ik}(z,\zeta)\Ab_{0,kj}^{lr}(\zeta)\d\zeta 
		\notag \\ 
		&\quad 	+ \sum_{k=1}^{n}\int_0^{1} P_{ik}(z,\zeta)\At_{0,kj}^{l}(\zeta)\d\zeta -\Ab_{0,ij}^{lr}(z)
		\label{eq:couplingCompBC2} 
	\end{align}
\end{subequations}
if the components of $\Af^{r}_i(z)$, $i=0,1$, are chosen to be zero for $i\leq j$ and 
\begin{subequations}\label{eq:newCouplingMatrices}
	\begin{flalign}
	& \underline{i>j:} \notag \\
		&\ \check{A}_{0,ij}^{r}(z) 
		 = 
		\lambda_j^{r}(0)\Q_{\zeta}(z,0) + \lambda_j^{r\prime}(0)\Q(z,0) + \lambda_j^{l}(0)\P_{\zeta}(z,0)  \notag \\ 
		&\ \quad + \lambda_j^{l\prime}(0) \P(z,0) - \sum_{k=1}^{n}\int_0^{z}Q_{ik}(z,\zeta)\Ab_{0,kj}^{lr}(\zeta)\d\zeta 
		\notag \\ 
		&\  \quad - \sum_{k=1}^{n}\int_0^{1} P_{ik}(z,\zeta)\At_{0,kj}^{l}(\zeta)\d\zeta +\Ab_{0,ij}^{lr}(z) \\
		&\ \check{A}_{1,ij}^{r}(z) 
		 =\ytn\Big( \lambda_j^{l}(0)\P(z,0) - \tfrac{1}{\ytn} \lambda_j^{r}(0)\Q(z,0) \notag \\ 
		&\ \quad + \sum_{k=1}^{n}\int_0^{z}Q_{ik}(z,\zeta) \Ab_{1,kj}^{lr}(\zeta) \d\zeta 
		\notag \\ 
		&\	\quad + \sum_{k=1}^{n}\int_0^{1}P_{ik}(z,\zeta) \At_{1,kj}^{l}(\zeta)\d\zeta - \Ab_{1,ij}^{lr}(z) \Big).
	\end{flalign}
\end{subequations}
\mrk{Note that this degree of freedom is needed to ensure the well-posedness of the kernel equations for the Volterra-kernel $Q(z,\zeta)$, just like it was for the backstepping kernel $K(z,\zeta)$ (see Section \ref{sec:canonKernelBackstepping}). 
Since both kernels of the second transformation are of dimension $n\times n$, the strictly lower triangular
well-posedness terms defined by \eqref{eq:newCouplingMatrices} are enough to ensure well-posedness of the kernel equations. This is a simplification compared to the first transformation with Figure \ref{fig:areas}.
}

While the BVP \eqref{eq:PCompPde}--\eqref{eq:P3} is already in the form which can be transformed into canonical coordinates, a further simplification is required for \eqref{eq:QCompPde}--\eqref{eq:Q4}.
In contrast to \cross{the kernel equations }\eqref{eq:kernelEquations}, the BVPs for $\Q(z,\zeta)$ contain no inhomogeneity except in the coupling BCs \eqref{eq:couplingCompBC1}, \eqref{eq:couplingCompBC2}. Thus, \eqref{eq:Q3} can be rewritten for $i=j$ to obtain
\begin{align}
	\Q'(z,z) = -\frac{\lambda_i^{r\prime}(z)}{2\lambda_i^{r}(z)}\Q(z,z)
\end{align}
with the solution $\Q(z,z) = 0$ following from \eqref{eq:Q4}.
Since $\lambda_i^{r} \neq \lambda_j^{r}$ holds for $i\neq j$, \eqref{eq:Q2} shows that $\Q(z,z)=0$ also holds in that case. Hence, the BVP \eqref{eq:QCompPde}--\eqref{eq:Q4} can be represented in simplified \textit{component form} as 
\begin{subequations}
	\begin{align}
		& \lambda_i^{r}(z)\Q_{zz}(z,\zeta) - (\lambda_j^{r}(\zeta)\Q(z,\zeta))_{\zeta\zeta} = 0 \label{eq:Q1'}\\
		& \Q(z,z) = 0 \label{eq:Q2'}\\
		& [\Q_z(z,z) = 0]_{\rightcond{i\neq j}} \label{eq:Q3'} 
	\end{align}%
\end{subequations}
complemented by the coupling BCs \eqref{eq:couplingCompBC1}, \eqref{eq:couplingCompBC2}.

\newcommand{\N}{\bm{N}}
\newcommand{\M}{\bm{M}}
\renewcommand{\H}{\pH}
\newcommand{\D}{\bm{D}}
\newcommand{\G}{\pG}
\begin{figure*}
	\centering
	\input{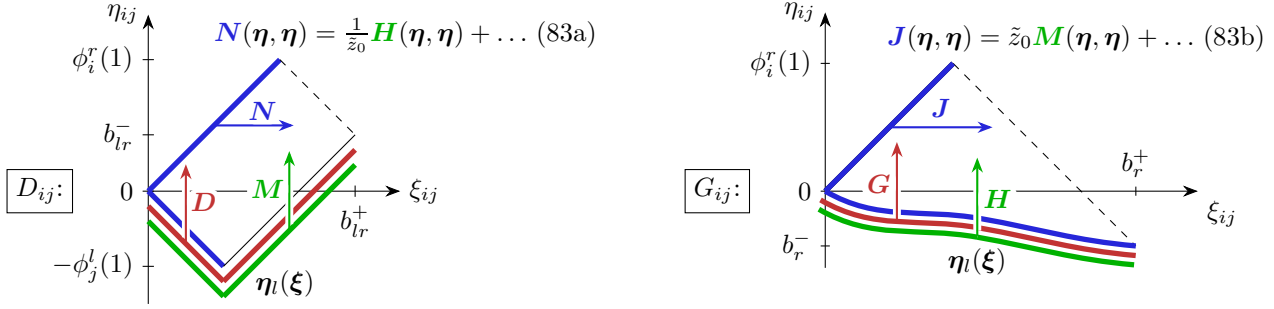}
	\caption{%
		The spatial domains of the canonical kernel equations \eqref{eq:canonD}, \eqref{eq:canonP} and \eqref{eq:K1K2TilTil} in the canonical coordinate systems $(\pxi,\peta)$ for $i\leq j$ with the utilized directions of integration. 
		The spatial domains are characterized by $b^{-}_{lr} = \phi_i^{r}(1) - \phi_j^{l}(1)$, $b^{+}_{lr}=\phi_i^{r}(1) + \phi_j^{l}(1)$, $b^- = \phi_i^{r}(1) - \phi_j^{r}(1)$, $b^+ = \phi_i^{r}(1)+\phi_j^{r}(1)$. In the left picture, the lower boundary $\peta_l(\pxi)$ is given by \eqref{eq:etaLFredholm}. In the right picture, it is determined the same way as in Figure \ref{fig:intDirections}. The thick coloured lines represent the BCs and the coloured arrows the respective directions of integration for the corresponding variables. Herein, the blue and red arrows symbolize the first and second integration, respectively. For the left elements, the BC for $\pJ(\peta,\peta)$ depends on $\H(\peta,\peta)$ for which the respective integral equation can be substituted, which is symbolized by the green arrow. On the other hand, $\pJ(\peta,\peta)$ depends on $\M(\peta,\peta)$ for which again the corresponding integral equation can be substituted.}
	\label{fig:intDirectionsDecoupling}
\end{figure*}

To convert the kernel equations into their canonical form, the \textit{canonical coordinates}
\newcommand{\xip}{\pxi}
\newcommand{\xipFun}{\pxi^{P}}
\newcommand{\etap}{\peta}
\newcommand{\etapFun}{\peta^{P}}
\newcommand{\xiq}{\pxi}
\newcommand{\xiqFun}{\pxi^{Q}}
\newcommand{\etaq}{\peta}
\newcommand{\etaqFun}{\peta^{Q}}
\begin{subequations}\label{eq:xieta(z,zeta)Decoupling}
	\begin{align}
			\xi_{ij}^{P}(z,\zeta) &=	\phi_i^{r}(z)+\phi_j^{l}(\zeta) \label{eq:xi(z,zeta)P} \\
			\eta_{ij}^{P}(z,\zeta) &= \phi_i^{r}(z) - \phi_j^{l}(\zeta) \label{eq:eta(z,zeta)P}	\\
			\xi_{ij}^{Q}(z,\zeta) &= \tfrac{1}{2}(1-\s)(\phi_i^{r}(1)+\phi_j^{r}(1)) + 
			\s(\phi_i^{r}(z)+\phi_j^{r}(\zeta)) \label{eq:xi(z,zeta)Q}\\
			\eta_{ij}^{Q}(z,\zeta) &= 
			-\tfrac{1}{2}(1-\s)(\phi_i^{r}(1)-\phi_j^{r}(1)) + \phi_i^{r}(z) - \phi_j^{r}(\zeta)	\label{eq:eta(z,zeta)Q}
	\end{align}
\end{subequations}
with 
\begin{subequations}
	\begin{align}
		\phi_i^{r}(z) = \int_0^{z}\frac{\d\zeta}{\sqrt{\lambda_i^{r}(\zeta)}}, \quad 
		\phi_i^{l}(z) = \int_0^{z}\frac{\d\zeta}{\sqrt{\lambda_i^{l}(\zeta)}},
	\end{align}
\end{subequations}
$i=1,\ldots,n$ and the corresponding inverses 
\newcommand{\zp}{z^{P}}
\newcommand{\zq}{z^{Q}}
\newcommand{\zetap}{\zeta^{P}}
\newcommand{\zetaq}{\zeta^{Q}}
\begin{subequations}\label{eq:zZetaInvP}
	\begin{align}
		\zp(\pxi,\peta) &=  (\phi_i^r)^{-1}(\tfrac{1}{2}(\pxi+\peta)) \\
		\zetap(\pxi,\peta) &= (\phi_j^l)^{-1}(\tfrac{1}{2}(\pxi-\peta))
\end{align}
\end{subequations}
as well as 
\begin{subequations}\label{eq:zZetaInvQ}
	\begin{align}
		\zq(\pxi,\peta) &=  (\phi_i^r)^{-1}(\tfrac{1}{2}(\s\pxi+\peta)+\tfrac{1}{2}(1-\s)\phi_i^{r}(1)) \\
		\zetaq(\pxi,\peta) &= (\phi_j^r)^{-1}(\tfrac{1}{2}(\s\pxi-\peta)+\tfrac{1}{2}(1-\s)\phi_j^{r}(1)),
	\end{align}
\end{subequations}
respectively,
are introduced along with the new kernel elements
\begin{subequations}\label{eq:DG}%
	\begin{flalign}%
		\D(\pxi,\peta) &= D_{ij}(\xipFun(z,\zeta),\etapFun(z,\zeta)) = \lambda_j^{l}(\zeta) \P(z,\zeta) \hspace{-5mm}&\\
		\G(\pxi,\peta) &= G_{ij}(\xiqFun(z,\zeta),\etaqFun(z,\zeta)) = \lambda_j^{r}(\zeta) \Q(z,\zeta).\hspace{-1cm}&
	\end{flalign}
\end{subequations}
With the change of coordinates \eqref{eq:xi(z,zeta)P}, \eqref{eq:eta(z,zeta)P}, the boundaries of the original Fredholm-domain of $P(z,\zeta)$ are mapped as $(z,0)\to (\peta,\peta)$, $(0,\zeta)\to(\peta,-\peta)$, $(z,1) \to (\pxi,\pxi-2\phi_j^{l}(1))$ and $(1,\zeta)\to (\pxi,2\phi_i^{r}(1)-\pxi)$. Hence the lower boundary 
\begin{align}\label{eq:etaLFredholm}
	\peta_l(\pxi) = \peta_l^{D}(\pxi) = \begin{cases}
		-\pxi, & \pxi \leq \phi_j^{l}(1) \\
		\pxi-2\phi_j^{l}(1), & \pxi > \phi_j^{l}(1)
	\end{cases}
\end{align}
consists of two parts in this case. Of course, the same is valid for the left boundary
\begin{align}
	\pxi_l(\peta) = \pxi_l^{D}(\peta) = \begin{cases}
		\peta, & \peta \geq 0 \\
		-\peta, & \peta < 0.
	\end{cases}
\end{align}%
The resulting spatial domains for $\D$ and $\G$ are depicted in Figure \ref{fig:intDirectionsDecoupling}.
Due to the fact that \eqref{eq:Q1'} only contains the diffusion coefficients of the right subsystem, the result $\lambda_i^{r} \geq \lambda_j^{r}$ automatically holds for $i\leq j$, which is why specifying the matrices $\Af_i^{r}(z)$, $i=0,1$, as strictly lower triangular according to \eqref{eq:newCouplingMatrices} leads to a solvable BVP for $\Q(z,\zeta)$.
To ensure the same shape of the spatial domain for all $i,j$ (see Figure \ref{fig:intDirections}), $s_{ij}$ according to \eqref{shdef} is introduced in the transformation \eqref{eq:xi(z,zeta)Q}, \eqref{eq:eta(z,zeta)Q} like it was applied in \eqref{eq:xieta(z,zeta)}.
In contrast to that, \eqref{eq:PCompPde} contains $\lambda_i^{l}$ and $\lambda_j^{r}$. In this case, $\lambda_i^{r} < \lambda_j^{l}$ for all $i,j=1,\ldots,n$ due to the assumed sorting of the diffusion coefficients. Hence the resulting spatial domains for all kernel elements $\D$, \ie\ the Fredholm part, already have a similar shape. Moreover, due to the types of BCs, it can be seen from the transformation into integral equations that the simple change of coordinates \eqref{eq:xi(z,zeta)Q}, \eqref{eq:eta(z,zeta)Q} is sufficient for the Fredholm kernel.

\mrk{Together, \eqref{eq:xieta(z,zeta)Decoupling}, \eqref{eq:DG} lead} to the \textit{canonical kernel equations}
\begin{subequations}\label{eq:canonD}
	\begin{flalign}
			&\D_{\xip\etap}(\xip,\etap)= 
			 \underbrace{-\tfrac{\ad^{lr}}{4} \D_{\xip}(\xip,\etap) + \tfrac{\as^{lr}}{4} \D_{\etap}(\xip,\etap)}_{
			\bm{\mA}_{\D}[\D_{\xi},\D_{\eta}](\pxi,\peta)
			} \label{eq:canonDPDE} &\\
			&\D(\etap,-\etap) = 0 \label{eq:DBC1} &\\
			&\D_{\etap}(\etap,-\etap) = 0  \label{eq:DBC2}&\\
			&\D_{\xip}(\etap,-\etap) = 0  \label{eq:DBC2'}&\\
			& \big[\D_{\xip}(\xip,\xip-2\phi_j^{l}(1)) = \D_{\etap}(\xip,\xip-2\phi_j^{l}(1))\big]_{\rightcond{\xip>\phi_j^{l}(1)}} \label{eq:DBC3} \hspace{-7mm}&
	\end{flalign}%
\end{subequations}%
and
\begin{subequations}\label{eq:canonP}
	\begin{align}
		&\G_{\xiq\etaq}(\xiq,\etaq)= 
		\underbrace{-\tfrac{\ad^{r}}{4} \G_{\xiq}(\xiq,\etaq) + \tfrac{\as^{r}}{4} \G_{\etap}(\xiq,\etaq)}_{
			\bm{\mA}_{\G}[\G_{\xi},\G_{\eta}](\pxi,\peta)
		} \label{eq:canonPPDE}\\
		& \G(\xiq,\etaq_l(\xiq)) = 0 \label{eq:GBC1}\\
		& \G_{\pxi}(\xiq,\etaq_l(\xiq)) = 0 \label{eq:GBC2} \\
		& [\G_{\peta}(\xiq_l(\etaq),\etaq) = 0]_{\rightcond{\etaq < 0}} \label{eq:GBC3}
	\end{align}%
\end{subequations}%
\mrk{with the coupling BCs}
\begin{subequations}\label{eq:canonCouplingFirst}
	\begin{flalign}
		&\underline{i\leq j:} \notag \\
		&\D(\etap,\etap) = \tfrac{1}{\ytn}\G(\etaq,\etaq) + \bm{\Ab}_1^{lr}(\zp(\etap,\etap)) \notag \\ 
		&	- \int_0^{\zp(\etap,\etap)}
			\frac{1}{\lambda_k^{r}(\zeta)}G_{ik}((\xi^{Q}_{ik},\eta^{Q}_{ik})(\zp(\etap,\etap),\zeta)) \Ab_{1,kj}^{lr}(\zeta)\d\zeta 
		\notag \\ 
		&- \int_0^{1}\!\frac{1}{\lambda_k^{l}(\zeta)}D_{ik}((\xi^{P}_{ik},\eta^{P}_{ik})(\zp(\etap,\etap),\zeta)) \At_{1,kj}^{l}(\zeta)\d\zeta \hspace{-5mm}& \label{eq:K2Til}\\[5pt]
		&\G_{\xip}(\etap,\etap) - \G_{\etap}(\etap,\etap) =
		-\ytn (\D_{\xiq}(\etaq,\etaq) - \D_{\etaq}(\etaq,\etaq)) \notag \\ 
		& + \int_0^{\zq(\etaq,\etaq)}
			\frac{\sqrt{\ljLow(0)}}{\lambda_k^{r}(\zeta)}G_{ik}((\xi^{Q}_{ik},\eta^{Q}_{ik})(\zq(\etaq,\etaq),\zeta))
			\Ab_{0,kj}^{lr}(\zeta) \d\zeta \notag \\ 
			& + \int_0^{1} \frac{\sqrt{\ljLow(0)}}{\lambda_k^{l}(\zeta)} D_{ik}((\xi^{P}_{ik},\eta^{P}_{ik})(\zp(\etaq,\etaq),\zeta)) \At_{0,kj}^{l}(\zeta)\d\zeta \notag \\ 
			& - \sqrt{\ljLow(0)}\bm{\Ab}_0^{lr}(\zp(\etap,\etap)),
			\label{eq:K1Til}
	\end{flalign}
\end{subequations}
where
\begin{subequations}
	\begin{align}
		\ad^{lr} &= \ad^{lr}(\pxi,\peta) = \tfrac{\lambda_j^{l\prime}(\zeta)}{2\sqrt{\lambda_j^{l}(\zeta)}}
		- \tfrac{\lambda_i^{r\prime}(z)}{\sqrt{\lambda_i^{r}(z)}} \\
		\as^{lr} &= \as^{lr}(\pxi,\peta) = \tfrac{\lambda_j^{l\prime}(\zeta)}{2\sqrt{\lambda_j^{l}(\zeta)}}
		+ \tfrac{\lambda_i^{r\prime}(z)}{\sqrt{\lambda_i^{r}(z)}}
	\end{align}
\end{subequations}%
and $(z,\zeta)$ are substituted by \eqref{eq:zZetaInvP} as well as%
\begin{subequations}
	\begin{align}
		\ad^{r} &= \ad^{r}(\pxi,\peta) = \tfrac{\lambda_j^{r\prime}(\zeta)}{2\sqrt{\ljLow(\zeta)}}
		- \tfrac{\lambda_i^{r\prime}(z)}{\lambda_i^{r}(z)} \\
		\as^{r} &= \as^{r}(\pxi,\peta) = \tfrac{\lambda_j^{r\prime}(\zeta)}{2\sqrt{\ljLow(\zeta)}}
		+ \tfrac{\lambda_i^{r\prime}(z)}{\lambda_i^{r}(z)},
	\end{align}
\end{subequations}
in which $(z,\zeta)$ are substituted by \eqref{eq:zZetaInvQ}.
\mrk{To simplify the notation, the \textit{summation convention} $\sum_{k=1}^{n} (c_k A_{ik}B_{kj}) \eqqcolon c_k A_{ik}B_{kj}$ is introduced in \eqref{eq:canonCouplingFirst}, which means that all expressions having $k$ as an index, are summed from $1$ to $n$.}
To be able to convert the canonical kernel equations into integral equations, the coupling BCs \eqref{eq:canonCouplingFirst} require an additional reformulation. To this end, \eqref{eq:K2Til} is differentiated \wrt\ $\peta$. The result is utilized in \eqref{eq:K1Til} to obtain the BCs%
\begin{subequations}\label{eq:K1K2TilTil}%
	\begin{align}%
		&\underline{i\leq j:} \notag \\
		&\D_{\peta}(\peta,\peta) = \tfrac{1}{\ytn} \G_{\pxi}(\peta,\peta) 
		- \tfrac{1}{2\ytn}\int_0^{z} \frac{\sqrt{\ljLow(0)}}{\lambda_k^{r}(\zeta)} 
			G_{ik} \Ab_{0,kj}^{lr}(\zeta) \d\zeta \notag \\ 
		& 
		- \tfrac{1}{2\ytn} \int_0^{1}\frac{\sqrt{\ljLow(0)}}{\lambda_k^{l}(\zeta)}
			D_{ik} \At_{0,kj}^{l}(\zeta) \d\zeta  \notag \\ 
		&
		- \int_0^{z}\frac{1}{2\lambda_k^{r}(\zeta)}(s_{ik} G_{ik,\xi} + G_{ik,\eta}) \Ab_{1,kj}^{lr}(\zeta) \d\zeta \notag \\ 
		& 
		- \int_0^{1}\frac{1}{2\lambda_k^{l}(\zeta)}(D_{ik,\xi} + D_{ik,\eta}) \At_{1,kj}^{l}(\zeta) \d\zeta \notag \\ 
		&
		- \frac{\sqrt{\ljLow(0)}}{2\ytn}\bm{\Ab}_0^{lr}(z)
		+ \frac{\sqrt{\lambda_i^{r}(z)}}{2}(\bm{\Ab}_1^{lr})'(z)
		\label{eq:canonCouplingBC1} \\
		&\G_{\peta}(\peta,\peta) = \ytn \D_{\pxi}(\peta,\peta) 
		- \int_0^{z} \frac{\sqrt{\ljLow(0)}}{2\lambda_k^{r}(\zeta)} G_{ik} \Ab_{0,kj}^{lr}(\zeta)\d\zeta \notag \\ 
		& 
		- \int_0^{1} \frac{\sqrt{\ljLow(0)}}{2\lambda_k^{l}(\zeta)} D_{ik} \At_{0,kj}^{l}(\zeta)\d\zeta \notag \\ 
		& 
		+ \int_0^{z} \frac{\ytn}{2\lambda_k^{r}(\zeta)} (s_{ik} G_{ik,\xi} + G_{ik,\eta}) \Ab_{1,kj}^{lr}(\zeta)\d\zeta \notag \\ 
		& 
		+ \int_0^{1} \frac{\ytn}{2\lambda_k^{l}(\zeta)} (D_{ik,\xi} + D_{ik,\eta}) \At_{1,kj}^{l} (\zeta)\d\zeta \notag \\ 
		&
		+ \frac{\sqrt{\ljLow(0)}}{2}\bm{\Ab}_0^{lr}(z) - \frac{\sqrt{\lambda_i^{r}(z)}}{2} (\bm{\Ab}_1^{lr})'(z).
		\label{eq:canonCouplingBC2}
	\end{align}
\end{subequations}
In \eqref{eq:canonCouplingBC1} $z=\zp(\peta,\peta)$, $G_{ik}$ and its derivatives have the argument $(\xi^{Q}_{ik},\eta^{Q}_{ik})(\zp(\etap,\etap),\zeta)$ and  $D_{ik}$ and its derivatives are evaluated at $(\xi^{P}_{ik},\eta^{P}_{ik})(\zp(\etap,\etap),\zeta)$. In \eqref{eq:couplingBC2} $z=\zq(\peta,\peta)$, $G_{ik}$ and its derivatives have the argument $(\xi^{Q}_{ik},\eta^{Q}_{ik})(\zq(\etap,\etap),\zeta)$, and  $D_{ik}$ and its derivatives are evaluated at $(\xi^{P}_{ik},\eta^{P}_{ik})(\zq(\etap,\etap),\zeta)$.

By mapping $(\bm{\Ab}_1^{lr})'(z)$ in \eqref{eq:At1ij} (see \eqref{eq:couplingMatrices}) to canonical coordinates, it can be shown that it is a piecewise continuous function.

To uniquely determine the kernels, the BCs are complemented by the \textit{artificial BCs}
\begin{subequations}\label{eq:artBCDecoupling}
	\begin{align}
			&\mathclap{\underline{i>j:}} \notag \\ 
		  &\D_{\peta}(\peta,\peta) = \bm{g}_{D}(\peta) \\
		  &\G_{\peta}(\peta,\peta) = \bm{g}_{G}(\peta)
	\end{align}
\end{subequations}
(see \eqref{eq:wellPosedBC}, \eqref{eq:wellPosedBC2}) with the degrees of freedom $\bm{g}_{D} = g_{D,ij} \in C[0, \phi_i^{r}(1)]$ and $\bm{g}_{G} = g_{G,ij} \in C[0,\phi_i^{r}(1)]$, $i,j=1,\ldots,n$. 

\subsubsection{Kernel integral equations}
Similar to Section \ref{sec:kernel integral equations}, the canonical kernel equations \eqref{eq:canonD}, \eqref{eq:canonP} with \eqref{eq:K1K2TilTil}\mrk{, \eqref{eq:artBCDecoupling}} are converted into integral equations by formally integrating the PDEs \eqref{eq:canonDPDE}, \eqref{eq:canonPPDE} \wrt\ $\pxi$ and $\peta$. \cross{In the present case,} \mrk{A} BC at $(\pxi,\peta) = (\peta,\peta)$ for the respective derivative \wrt\ $\peta$ is available for both $\D$ and $\G$ by \eqref{eq:K1K2TilTil}. Together with \eqref{eq:DBC2}, \eqref{eq:GBC3} \mrk{and \eqref{eq:canonCouplingFirst}}, $\D_{\peta}(\pxi_l(\peta),\peta)$ and $\G_{\peta}(\pxi_l(\peta),\peta)$ are known at the whole left boundary $\pxi_l(\peta)$ of the respective spatial domain (see Figure \ref{fig:intDirectionsDecoupling}). 

To get the BC for $\D$ and $\G$ on the lower boundary, note that $\D'(\pxi,\pxi-2\phi_j^{l}(1))$ = $\D_{\pxi}(\pxi,\pxi-2\phi_j^{l}(1)) + \D_{\peta}(\pxi,\pxi-2\phi_j^{l}(1))$. Inserting \eqref{eq:DBC3} and integrating \wrt\ $\pxi$ then yields
\begin{align}\label{eq:Dxil}
	&\big[\D(\pxi,\pxi-2\phi_j^{l}(1)) = \underbrace{\D(\phi_j^{l}(1),-\phi_j^{l}(1))}_{\stackrel{\eqref{eq:DBC2}}{=}0}
	\notag \\ 
	&\quad   	+ \int_{\phi_j^{l}(1)}^{\pxi} 2\D_{\peta} (\oxi,\oxi-2\phi_j^{l}(1)) \d\oxi\, \big]_{\rightcond{\pxi > \phi_j^{l}(1)}}.
\end{align}
Together with \eqref{eq:DBC2} and \eqref{eq:GBC1}, BCs $\D(\pxi,\peta_l(\pxi))$ and $\G(\pxi,\peta_l(\pxi))$ are available at the whole lower boundary $\peta_l(\pxi)$ of the respective spatial domain (see Figure \ref{fig:intDirectionsDecoupling}). 

Furthermore, note that even for the respective derivatives $\D_{\pxi}$ and $\G_{\pxi}$, a BC on the whole lower boundary of the domains is determined by \eqref{eq:DBC2'}, \eqref{eq:DBC3} and \eqref{eq:GBC2}.

Hence, introducing the variables 
\begin{subequations}
	\begin{align}
		\N(\pxi,\peta) &= N_{ij}(\pxi,\peta) \coloneqq \D_{\peta}(\pxi,\peta) \\
		\J(\pxi,\peta) &= J_{ij}(\pxi,\peta) \coloneqq \G_{\peta}(\pxi,\peta)
	\end{align}
\end{subequations}
	and 
\begin{subequations}
	\begin{align}
		\M(\pxi,\peta) &= M_{ij}(\pxi,\peta) \coloneqq \D_{\pxi}(\pxi,\peta) \\
		\pH(\pxi,\peta) &= G_{ij}(\pxi,\peta) \coloneqq \G_{\pxi}(\pxi,\peta),
	\end{align}
\end{subequations}
and inserting the BCs \eqref{eq:DBC1}--\eqref{eq:DBC3}, \eqref{eq:GBC1}--\eqref{eq:GBC3}, \eqref{eq:Dxil}, \eqref{eq:K1K2TilTil} and \eqref{eq:artBCDecoupling}, the \textit{kernel integral equations} finally read
\begin{subequations}\label{eq:inteqDecouplingFormal}
	\begin{align}
		\N &= \N_0 + \F_{N}[N,D,M,J,G,H] \label{eq:inteqNFormal} \\
		\D &= \D_0 + \F_{D}[\N] \\
		\M &= \M_0 + \F_{M}[N,D,M,J,G,H] \\
		\J &= \pJ_0 + \F_J[N,D,M,J,G,H] \\
		\G &= \F_{G}[\J] \\
		\H &= \F_H[\J,\H]
	\end{align}
\end{subequations}
with 
\begin{subequations}\label{eq:startValsDecoupling}
	\begin{align}
		&\N_0 = [\bm{g}_{D}(\peta)]_{\rightcond{i&>j \\ \peta& \geq 0}}  \\ 
		&  \quad
		+ \big[\frac{\sqrt{\ljLow(0)}}{2\ytn} \bm{\Ab}_0^{lr}(z) + 
		\frac{\sqrt{\li(z)}}{2}(\bm{\Ab}_1^{lr})'(z)
		\big]_{\rightcond{i &\leq j \\ z &= \zp(\peta,\peta) \\ \peta &\geq 0}} \notag\\
		& \D_0 = \Big[
			\int_{\phi_j^{l}(1)}^{\pxi} 2 \N_0(\peta_l(\oxi))\d\oxi
		\Big]_{\rightcond{\pxi>\phi_j^{l}(1)}} \\
		& \M_0 = [\N_0(\peta_l(\pxi))]_{\rightcond{\pxi&>\phi_j^{l}(1)}} \\
		& \J_0 = [\bm{g}_G(\peta)]_{\rightcond{i&>j \\ \peta& \geq 0}} +\Big[\ytn \M_0(\peta)   \\* 
		& \quad + \frac{\sqrt{\ljLow(0)}}{2}\bm{\Ab}_{0}^{lr}(z) - \frac{\ytn \sqrt{\li(z)}}{2} (\bm{\Ab}_1^{lr})'(z)\Big]_{\rightcond{i &\leq j \\ z &= \zq(\peta,\peta) \\ \peta &\geq 0}}\notag
	\end{align}
\end{subequations}
and
\begin{subequations}\label{eq:intOpsDecoupling}
	\begin{align}
		&\F_N[N,D,M,J,G,H]  =
		\big[
			\frac{1}{\ytn}\int_{\peta_l(\peta)}^{\peta}\bm{\mA}_G[\H,\J](\peta,\oeta)\d\oeta \notag \\ 
			& \quad
			- \frac{1}{2\ytn}\int_0^{z} \frac{\sqrt{\ljLow(0)}}{\lambda_k^{r}(\zeta)} 
			G_{ik} \Ab_{0,kj}^{lr}(\zeta) \d\zeta \notag \\ 
			& \quad
			- \frac{1}{2\ytn} \int_0^{1}\frac{\sqrt{\ljLow(0)}}{\lambda_k^{l}(\zeta)}
			D_{ik} \At_{0,kj}^{l}(\zeta) \d\zeta  \notag \\ 
			&\quad
			- \int_0^{z}\frac{1}{2\lambda_k^{r}(\zeta)}(s_{ik} H_{ik} + J_{ik}) \Ab_{1,kj}^{lr}(\zeta) \d\zeta \notag \\ 
			& \quad
			- \int_0^{1}\frac{1}{2\lambda_k^{l}(\zeta)}(M_{ik} + N_{ik}) \At_{1,kj}^{l}(\zeta) \d\zeta
		\big]_{\rightcond{i &\leq j \\ z&=\zp(\peta,\peta) \\ \peta&\geq 0}} \notag \\ 
		& \quad
			+ \int_{\pxi_l(\peta)}^{\pxi}\bm{\mA}_{D}[\M,\N](\oxi,\peta)\d\oxi \\
		&\F_D[\N] = \int_{\peta_l(\pxi)}^{\peta}\N(\pxi,\oeta)\d\oeta \notag \\ 
		& \quad + \Big[
			\int_{\phi_j^{l}(1)}^{\pxi} 2 \F_{N}(\oxi,\peta_l(\oxi))\d\oxi
		\Big]_{\rightcond{\pxi>\phi_j^{l}(1)}} \\
		& \F_M[N,D,M,J,G,H] = [\F_N(\pxi,\peta_l(\pxi))]_{\rightcond{\pxi&> \phi_j^{l}(1)}} \notag \\ 
		& \quad + \int_{\peta_l(\pxi)}^{\peta}\bm{\mA}_D[\M,\N](\pxi,\oeta)\d\oeta \\
		&\F_J[N,D,M,J,G,H]  =
		\big[
		\ytn \F_{M}(\peta,\peta) \notag \\ 
		& \quad
		- \frac{1}{2}\int_0^{z} \frac{\sqrt{\ljLow(0)}}{\lambda_k^{r}(\zeta)} 
		G_{ik} \Ab_{0,kj}^{lr}(\zeta) \d\zeta \notag \\ 
		& \quad
		- \frac{1}{2} \int_0^{1}\frac{\sqrt{\ljLow(0)}}{\lambda_k^{l}(\zeta)}
		D_{ik} \At_{0,kj}^{l}(\zeta) \d\zeta  \notag \\ 
		&\quad
		+ \int_0^{z}\frac{\ytn}{2\lambda_k^{r}(\zeta)}(s_{ik} H_{ik} + J_{ik}) \Ab_{1,kj}^{lr}(\zeta) \d\zeta \notag \\ 
		& \quad
		+ \int_0^{1}\frac{\ytn}{2\lambda_k^{l}(\zeta)}(M_{ik} + N_{ik}) \At_{1,kj}^{l}(\zeta) \d\zeta
		\big]_{\rightcond{i &\leq j \\ z&=\zq(\peta,\peta) \\ \peta&\geq 0}} \notag \\ 
		& \quad
		+ \int_{\pxi_l(\peta)}^{\pxi}\bm{\mA}_{G}[\H,\J](\oxi,\peta)\d\oxi \\
		&\F_G[\J] = \int_{\peta_l(\pxi)}^{\peta} \J(\pxi,\oeta)\d\oeta \\
		& \F_H[\J,\H] = \int_{\peta_l(\pxi)}^{\peta} \bm{\mA}_G[\H,\J](\pxi,\oeta)\d\oeta.
	\end{align}%
\end{subequations}%
Therein, the integral equation \eqref{eq:inteqNFormal} has been utilized in \eqref{eq:DBC3} and \eqref{eq:Dxil}. Moreover, the arguments of $D_{ik}$ and its derivatives as well as $G_{ik}$ and its derivatives are the same as in \eqref{eq:K1K2TilTil}.

\subsection{Successive approximation}
The integral equations \eqref{eq:intEqs} with \eqref{eq:startVals} and \eqref{eq:intOps} as well as \eqref{eq:inteqDecouplingFormal} with \eqref{eq:startValsDecoupling} and \eqref{eq:intOpsDecoupling} are of a similar form as in \cite{Deu17} and can now be solved by the \emph{method of successive approximations}, \ie\ by applying fixed-point iteration. Thereby, the solutions may be represented by
\begin{flalign}\label{eq:fixedPointIteration}
	\pG &= \sum_{l=0}^{\infty}\Delta \pG^l, \quad \pH = \sum_{l=0}^{\infty} \Delta \pH^l, \quad 
	\pJ = \sum_{l=0}^{\infty} \Delta \pJ^l &\\
	\D &= \sum_{l=0}^{\infty}\Delta \D^l, \quad \M = \sum_{l=0}^{\infty} \Delta \M^l, \quad 
	\N = \sum_{l=0}^{\infty} \Delta \N^l \hspace{-5mm}&
\end{flalign}
with $\Delta G^l = [\Delta G^l_{ij}]$, $\Delta H^l = [\Delta H^l_{ij}]$ and $\Delta J^l = [\Delta J^l_{ij}]$, $i,j=1,\ldots,2n$ in the case of \eqref{eq:intEqs} and $i,j=1,\ldots,n$ in the case of \eqref{eq:inteqDecouplingFormal} as well as $\Delta D^l = [\Delta D^l_{ij}]$, $\Delta M^l = [\Delta M^l_{ij}]$ and $\Delta N^l = [\Delta N^l_{ij}]$, $i,j=1,\ldots,n$,
which are calculated by the update law%
\begin{subequations}\label{eq:updateLaw}
	\begin{flalign}
		\Delta\pG^{l+1} &= {\bm{F}}_G[\Delta \pH^l,\Delta J^l],\ \Delta \pG^0 = \pG_0(\pxi,\peta) \label{eq:updateLawG}&\\
		 \Delta \pH^{l+1} &= {\bm{F}}_{H}[\Delta G^l,\Delta \pH^l,\Delta J^l],\ \Delta \pH^0 = \pH_0(\pxi,\peta) \label{eq:updateLawH} \hspace{-10mm}&\\
		 \Delta \pJ^{l+1} &= {\bm{F}}_{J}[\Delta G^l,\Delta H^l,\Delta J^l],\ \Delta \pJ^0 = \pJ_0(\pxi,\peta) \hspace{-3mm}& \label{eq:updateLawJ}
	\end{flalign}
\end{subequations}
in the case of \eqref{eq:intEqs} and 
similar for \eqref{eq:inteqDecouplingFormal}.
Of course, the integral equations \eqref{eq:intEqs}--\eqref{eq:intOps} and \eqref{eq:inteqDecouplingFormal}--\eqref{eq:intOpsDecoupling} differ from the form in \cite{Deu17} due to the coupling BCs \eqref{eq:kernelFoldBC1}, \eqref{eq:kernelFoldBC2} and \eqref{eq:K1K2TilTil}, respectively. However, the types of the appearing terms in the integral operators \eqref{eq:intOps}\mrk{, \eqref{eq:intOpsDecoupling}} are the same except for the terms with the Fredholm integrals and the boundary integral $\int_{\phi_j^{l}(1)}^{\pxi}\ldots\d\oxi$. \cross{However,}Noting that inside the Fredholm integrals $z = z(\peta,\peta)$ always holds allows to determine integral estimates for both terms with a similar reasoning as in \cite{Deu17}. Hence, with the same approach as in the latter reference,
absolute and uniform convergence of the series \eqref{eq:fixedPointIteration} follows by proving
	\begin{align}\label{eq:grow}
		|\Delta\pX^l(\pxi(z,\zeta),\peta(z,\zeta))| \leq \frac{M^{l+1}}{l!}(z-\gamma \zeta)^l 
	\end{align}
for each $\Delta \pX \in \{ \Delta \pG, \Delta \pH, \Delta \pJ, \Delta\D, \Delta\M, \Delta\N \}$, some $M>0$ and
\begin{align}\label{eq:gamma}
	\gamma \in 
		\left(\max\limits_{i > j}
		\sqrt{\tfrac{\lambda_i({z}_{\underline{\Delta}})}{\lambda_j({z}_{\underline{\Delta}})}},\ 1\right),
\end{align}
where
$
	{z}_{\underline{\Delta}} = \operatorname{argmin}(|\lambda_i(z)-  \lambda_j(z)|)
$
is the point of minimal difference between the diffusion coefficients $\lambda_i$ and $\lambda_j$ (see \cite{Deu17}). Note that $\max_{i > j}
\sqrt{{\lambda_i^{r}({z}_{\underline{\Delta}})}/{\ljLow({z}_{\underline{\Delta}})}} \leq \max_{i > j}
\sqrt{{\lambda_i({z}_{\underline{\Delta}})}/{\lambda_j({z}_{\underline{\Delta}})}}$ holds so that \eqref{eq:gamma} is sufficient for the convergence proof of \eqref{eq:inteqDecouplingFormal}.
Therefore, \eqref{eq:fixedPointIteration} provides the piecewise continuous solution of the integral equations \eqref{eq:intEqs} and \eqref{eq:inteqDecouplingFormal}, which proves Theorems \ref{thm:kernelEquations} and \ref{thm:kernelEquationsDecoupling}. 
%

\section{Example}
\label{sec:simulation}
\comment{Wenn ich Fig.4 auf die naechste Seite schiebe, dann wandert die Diskussion auf diese hier. Sieht dann etwa so aus wie jetzt.}%
Consider a system \eqref{eq:origSys} consisting of two coupled PDEs with the parameters
\renewcommand{\arraystretch}{1.2}
\begin{subequations}
\begin{align}
	\check{\Lambda}(\y) &= \begin{bmatrix}
		z^{2}+2 & 0 \\ 0 & \e^{-z}+1/2
	\end{bmatrix} \\
	\check{A}(\y) &= \begin{bmatrix}
		1 & 1+\y \\ \frac{1}{2}+\y & 1
	\end{bmatrix},
\end{align} 
\end{subequations}
subject to Neumann BCs, \ie\ $B_0=B_1 = 0$ in \eqref{eq:origBC1} and \eqref{eq:origBC2}, which is open-loop unstable. For the bilateral controller design, a numerical analysis of the possible folding points to ensure non-intersection diffusion coefficients yields the intervals $\yn \in I_i$, $I_1 = (0,0.412)$, $I_2 = (0.424,0.488)$, $I_3 = (0.568,0.596)$ and $I_4 = (0.604, 1)$.
The folding point is set to $\yn = 0.325$, leading to descending diffusion coefficients $\lambda_1(z) > \ldots > \lambda_{4}(z)$ in the folded system.
After specifying $\mu =10$, the degrees of freedom $\bm{g}_f(\peta)$ in the well-posed BCs \eqref{eq:wellPosedBC}, \eqref{eq:wellPosedBC2} and $\bm{g}_D(\peta)$, $\bm{g}_G(\peta)$ in  \eqref{eq:artBCDecoupling} are set to zero for simplicity, which determines the kernels. The successive approximations \eqref{eq:fixedPointIteration} are implemented in \textsc{Matlab} and truncated after the maximum deviation $\max_{z,\zeta,i,j} |\Delta X_{ij}^l| \leq 10^{-3}$, $\Delta X \in \{\Delta G, \Delta H, \Delta J, \Delta D, \Delta M, \Delta N \}$, which lead to 9 iterations for $K(z,\zeta)$ and 8 iterations for $P(z,\zeta)$ and $Q(z,\zeta)$. The kernels are discretized by 51 points in each direction $z$, $\zeta$. 
In order to limit the number of resulting nodes,
the grids $(\pxi,\peta)$ in the canonical coordinates are resampled \mrk{to 100 nodes in $\pxi$-direction and the same node distance for $\peta$.}

\begin{figure}
	\centering
	\pgfplotsset{
		every axis post/.append style={
			width=0.4\linewidth,
			legend style={
				font=\footnotesize,
			},
			every axis x label/.append style={yshift=2mm},
			yticklabel style={inner sep=0pt,xshift=-4pt},
			xticklabel style = {inner xsep = 0pt},
			title style={inner sep=0pt},
		},%
		every axis plot/.append style={
			line width=1.5pt
		}
	}%
	\newcommand{\mymycmd}{-10t}%
	\colorlet{green}{greenCol}
	\colorlet{red}{redCol}
	\colorlet{blue}{blueCol}
%
%
\begin{tikzpicture}

\begin{axis}[%
width=0.333\linewidth,
height=0.35\linewidth,
at={(0\linewidth,0\linewidth)},
scale only axis,
xmin=0,
xmax=1,
xlabel style={font=\color{white!15!black}},
xlabel={$t$},
ymin=0,
ymax=1,
axis background/.style={fill=white},
title style={font=\bfseries},
title={$\|w(t)\|_{L_2}/\|w(0)\|_{L_2}$},
legend style={legend cell align=left, align=left, draw=white!15!black},
every axis title/.append style={font=\footnotesize},every axis label/.append style={font=\footnotesize},every tick label/.append style={font=\scriptsize},every x tick label/.style={},every y tick label/.style={},every z tick label/.style={},yticklabel style={/pgf/number format/fixed,/pgf/number format/precision=2},
]
\addplot [color=blue]
  table[]{plots/+data/1_Norms-1.tsv};
\addlegendentry{unilateral}

\addplot [color=green]
  table[]{plots/+data/1_Norms-2.tsv};
\addlegendentry{bilateral}

\addplot [color=red, dashed, line width=1.5pt]
  table[]{plots/+data/1_Norms-3.tsv};
\addlegendentry{$2.26\e^{\mymycmd}$}

\end{axis}
\end{tikzpicture}%
%
	\hfill%
	{%
		\pgfplotsset{%
			every axis post/.append style={%
				execute at end axis={
					\legend{}
					\node[red, anchor=west,font=\small] at (axis cs:0.18,2){$u_0(t)$ (bilateral)};
					\node[green, anchor=west,font=\small] at (axis cs:0.18,-2){$u_1(t)$ (bilateral)};
					\node[blue, anchor=west,font=\small] at (axis cs:0.025,-9){$u_1(t)$ (unilateral)};
				},
				xmax = 0.6,
				ymax = 4
			}%
		}%
%
%
\begin{tikzpicture}

\begin{axis}[%
width=0.333\linewidth,
height=0.35\linewidth,
at={(0\linewidth,0\linewidth)},
scale only axis,
xmin=0,
xmax=1,
xlabel style={font=\color{white!15!black}},
xlabel={$t$},
ymin=-16,
ymax=2,
axis background/.style={fill=white},
title style={font=\bfseries},
title={$u_0(t), u_1(t)$},
legend style={legend cell align=left, align=left, draw=white!15!black},
every axis title/.append style={font=\footnotesize},every axis label/.append style={font=\footnotesize},every tick label/.append style={font=\scriptsize},every x tick label/.style={},every y tick label/.style={},every z tick label/.style={},yticklabel style={/pgf/number format/fixed,/pgf/number format/precision=2},
]
\addplot [color=blue]
table[]{plots/+data/2_Inputs-1.tsv};
\addlegendentry{data1}

\addplot [color=blue, dashed]
table[]{plots/+data/2_Inputs-2.tsv};
\addlegendentry{data2}

\addplot [color=red]
table[]{plots/+data/2_Inputs-3.tsv};
\addlegendentry{data3}

\addplot [color=red, dashed]
table[]{plots/+data/2_Inputs-4.tsv};
\addlegendentry{data4}

\addplot [color=green]
table[]{plots/+data/2_Inputs-5.tsv};
\addlegendentry{data5}

\addplot [color=green, dashed]
table[]{plots/+data/2_Inputs-6.tsv};
\addlegendentry{data6}

\end{axis}
\end{tikzpicture}%
%
	}%
	\caption{Comparison of the bilateral with an unilateral backstepping controller. The left picture shows the $L_2$-norms of the closed-loop system in the bilateral and unilateral cases. The corresponding control efforts are depicted in the right picture, where the solid lines belong to the first element of the input $u_i(t)$, $i=0,1$, and the dashed lines to their second element.}%
	\label{fig:decay}%
\end{figure}
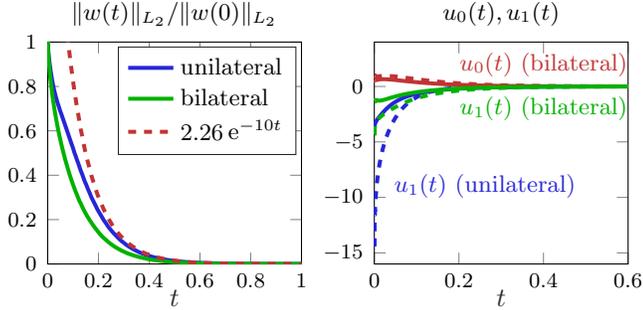
To evaluate the effect of the bilateral control, the designed state feedback controller \eqref{eq:controlLaw} is compared with an unilateral backstepping controller according to \cite{Deu17}, \ie\ $u_0(t) = 0$.
The left picture of Figure \ref{fig:decay} shows the resulting weighted $L_2$-norms $\|h\|_{L_2} = (\int_0^1\|\check{\Lambda}^{-\frac{1}{2}}(\y)h(\y)\|^2_{\mathbb{C}^n}\d \y)^{1/2}$ of the state $w(\y,t)$ 
with the 
IC $w_i(\y,0) = \frac{3}{4}\sin(\pi \y + 2\pi) + \frac{1}{4}\cos(3\pi \y + \frac{\pi}{2})$, $i=1,2$,
in both cases. It can be seen that the bilateral controller is able to exponentially stabilize the system with the desired decay rate. Moreover, transforming the resulting state profiles into the target coordinates by applying \eqref{eq:bsTrafo} and \eqref{eq:decouplingTrafo} and comparing with a simulation of the target system \eqref{eq:unfoldedTargetSystem} shows that the desired behaviour is achieved with a maximum deviation of $3.7\cdot 10^{-3}$, occurring solely due to numerics.


The right picture of Figure \ref{fig:decay} shows the corresponding control efforts in both cases. Due to the active left input $u_0(t)$ in the bilateral case, the required control effort $u_1(t)$ of the right input to achieve the desired decay rate is significantly reduced. 

\begin{figure}
	\colorlet{green}{greenCol}
	\colorlet{red}{redCol}
	\colorlet{blue}{blueCol}
	{
		\pgfplotsset{
			surf/.append style={shader = faceted,faceted color=black, fill=white},
			every axis/.append style={
				every axis y label/.append style={yshift=3mm, xshift=0mm},
				every axis x label/.append style={xshift=0mm,yshift=2mm},
				every axis z label/.append style={yshift=-2mm,xshift=0mm},
				extra z ticks={0},
				clip=false,
				xtick distance=0.5,
			},
			every axis post/.append style={
				width=0.38\linewidth,
				height=0.38\linewidth,
				xmajorgrids=false,
				ymajorgrids=false,
				zmajorgrids=false,
				tick align=inside,
				every axis y label/.append style={color=black},
				every axis x label/.append style={color=black},
				every axis z label/.append style={color=black},
				view  = {45}{25},
				zlabel={},
				zticklabel style={inner sep=0pt,xshift=-3pt},
				yticklabel style={inner xsep=0pt},
				ylabel = {$\y$},
				extra tick style={grid=major, ticklabel style={anchor=west, yshift=-1mm,xshift=1mm}},
				extra z ticks ={}, 
				title style={inner sep=0pt},
			},
			execute at end axis={
				\legend{}
			}
		}
%
%
\begin{tikzpicture}

\begin{axis}[%
width=0.262\linewidth,
height=0.35\linewidth,
at={(-0.012\linewidth,0\linewidth)},
scale only axis,
xmin=0,
xmax=1,
tick align=outside,
xlabel style={font=\color{white!15!black}},
xlabel={$t$},
ymin=0,
ymax=1,
ylabel style={font=\color{white!15!black}},
ylabel={$z$},
zmin=-0.2,
zmax=1,
zlabel style={font=\color{white!15!black}},
zlabel={$x_1(z,t)$},
view={-37.5}{30},
axis background/.style={fill=white},
extra y ticks = {0.16},
extra y tick labels = {$\textcolor{redCol}{\yn}$},
xmajorgrids,
ymajorgrids,
zmajorgrids,
legend style={at={(1.03,1)}, anchor=north west, legend cell align=left, align=left, draw=white!15!black},
every axis title/.append style={font=\footnotesize},every axis label/.append style={font=\footnotesize},every tick label/.append style={font=\scriptsize},every x tick label/.style={},every y tick label/.style={},every z tick label/.style={},yticklabel style={/pgf/number format/fixed,/pgf/number format/precision=2},
title = {$w_1(\y,t)$, $\yn = 0.16$},
]

\addplot3[%
surf,
shader=flat corner, draw=black, z buffer=none, colormap={mymap}{[1pt] rgb(0pt)=(0.2422,0.1504,0.6603); rgb(1pt)=(0.25039,0.164995,0.707614); rgb(2pt)=(0.257771,0.181781,0.751138); rgb(3pt)=(0.264729,0.197757,0.795214); rgb(4pt)=(0.270648,0.214676,0.836371); rgb(5pt)=(0.275114,0.234238,0.870986); rgb(6pt)=(0.2783,0.255871,0.899071); rgb(7pt)=(0.280333,0.278233,0.9221); rgb(8pt)=(0.281338,0.300595,0.941376); rgb(9pt)=(0.281014,0.322757,0.957886); rgb(10pt)=(0.279467,0.344671,0.971676); rgb(11pt)=(0.275971,0.366681,0.982905); rgb(12pt)=(0.269914,0.3892,0.9906); rgb(13pt)=(0.260243,0.412329,0.995157); rgb(14pt)=(0.244033,0.435833,0.998833); rgb(15pt)=(0.220643,0.460257,0.997286); rgb(16pt)=(0.196333,0.484719,0.989152); rgb(17pt)=(0.183405,0.507371,0.979795); rgb(18pt)=(0.178643,0.528857,0.968157); rgb(19pt)=(0.176438,0.549905,0.952019); rgb(20pt)=(0.168743,0.570262,0.935871); rgb(21pt)=(0.154,0.5902,0.9218); rgb(22pt)=(0.146029,0.609119,0.907857); rgb(23pt)=(0.138024,0.627629,0.89729); rgb(24pt)=(0.124814,0.645929,0.888343); rgb(25pt)=(0.111252,0.6635,0.876314); rgb(26pt)=(0.0952095,0.679829,0.859781); rgb(27pt)=(0.0688714,0.694771,0.839357); rgb(28pt)=(0.0296667,0.708167,0.816333); rgb(29pt)=(0.00357143,0.720267,0.7917); rgb(30pt)=(0.00665714,0.731214,0.766014); rgb(31pt)=(0.0433286,0.741095,0.73941); rgb(32pt)=(0.0963952,0.75,0.712038); rgb(33pt)=(0.140771,0.7584,0.684157); rgb(34pt)=(0.1717,0.766962,0.655443); rgb(35pt)=(0.193767,0.775767,0.6251); rgb(36pt)=(0.216086,0.7843,0.5923); rgb(37pt)=(0.246957,0.791795,0.556743); rgb(38pt)=(0.290614,0.79729,0.518829); rgb(39pt)=(0.340643,0.8008,0.478857); rgb(40pt)=(0.3909,0.802871,0.435448); rgb(41pt)=(0.445629,0.802419,0.390919); rgb(42pt)=(0.5044,0.7993,0.348); rgb(43pt)=(0.561562,0.794233,0.304481); rgb(44pt)=(0.617395,0.787619,0.261238); rgb(45pt)=(0.671986,0.779271,0.2227); rgb(46pt)=(0.7242,0.769843,0.191029); rgb(47pt)=(0.773833,0.759805,0.16461); rgb(48pt)=(0.820314,0.749814,0.153529); rgb(49pt)=(0.863433,0.7406,0.159633); rgb(50pt)=(0.903543,0.733029,0.177414); rgb(51pt)=(0.939257,0.728786,0.209957); rgb(52pt)=(0.972757,0.729771,0.239443); rgb(53pt)=(0.995648,0.743371,0.237148); rgb(54pt)=(0.996986,0.765857,0.219943); rgb(55pt)=(0.995205,0.789252,0.202762); rgb(56pt)=(0.9892,0.813567,0.188533); rgb(57pt)=(0.978629,0.838629,0.176557); rgb(58pt)=(0.967648,0.8639,0.16429); rgb(59pt)=(0.96101,0.889019,0.153676); rgb(60pt)=(0.959671,0.913457,0.142257); rgb(61pt)=(0.962795,0.937338,0.12651); rgb(62pt)=(0.969114,0.960629,0.106362); rgb(63pt)=(0.9769,0.9839,0.0805)}, mesh/rows=21]
table[point meta=\thisrow{c}] {%
	plots/+data/6_x_1-1.tsv};

\addplot3 [color=red, line width=2.0pt]
table[] {plots/+data/6_x_1-2.tsv};
\addlegendentry{data2}

\end{axis}
\end{tikzpicture}%
\hfill
%
%
\begin{tikzpicture}

\begin{axis}[%
width=0.262\linewidth,
height=0.35\linewidth,
at={(-0.012\linewidth,0\linewidth)},
scale only axis,
xmin=0,
xmax=1,
tick align=outside,
xlabel style={font=\color{white!15!black}},
xlabel={$t$},
ymin=0,
ymax=1,
ylabel style={font=\color{white!15!black}},
ylabel={$z$},
zmin=-0.2,
zmax=1,
zlabel style={font=\color{white!15!black}},
zlabel={$x_1(z,t)$},
view={-37.5}{30},
axis background/.style={fill=white},
xmajorgrids,
ymajorgrids,
zmajorgrids,
extra y ticks = {0.66},
extra y tick labels = {$\textcolor{redCol}{\yn}$},
legend style={at={(1.03,1)}, anchor=north west, legend cell align=left, align=left, draw=white!15!black},
every axis title/.append style={font=\footnotesize},every axis label/.append style={font=\footnotesize},every tick label/.append style={font=\scriptsize},every x tick label/.style={},every y tick label/.style={},every z tick label/.style={},yticklabel style={/pgf/number format/fixed,/pgf/number format/precision=2},
title = {$w_1(\y,t)$, $\yn = 0.66$},
]

\addplot3[%
surf,
shader=flat corner, draw=black, z buffer=none, colormap={mymap}{[1pt] rgb(0pt)=(0.2422,0.1504,0.6603); rgb(1pt)=(0.25039,0.164995,0.707614); rgb(2pt)=(0.257771,0.181781,0.751138); rgb(3pt)=(0.264729,0.197757,0.795214); rgb(4pt)=(0.270648,0.214676,0.836371); rgb(5pt)=(0.275114,0.234238,0.870986); rgb(6pt)=(0.2783,0.255871,0.899071); rgb(7pt)=(0.280333,0.278233,0.9221); rgb(8pt)=(0.281338,0.300595,0.941376); rgb(9pt)=(0.281014,0.322757,0.957886); rgb(10pt)=(0.279467,0.344671,0.971676); rgb(11pt)=(0.275971,0.366681,0.982905); rgb(12pt)=(0.269914,0.3892,0.9906); rgb(13pt)=(0.260243,0.412329,0.995157); rgb(14pt)=(0.244033,0.435833,0.998833); rgb(15pt)=(0.220643,0.460257,0.997286); rgb(16pt)=(0.196333,0.484719,0.989152); rgb(17pt)=(0.183405,0.507371,0.979795); rgb(18pt)=(0.178643,0.528857,0.968157); rgb(19pt)=(0.176438,0.549905,0.952019); rgb(20pt)=(0.168743,0.570262,0.935871); rgb(21pt)=(0.154,0.5902,0.9218); rgb(22pt)=(0.146029,0.609119,0.907857); rgb(23pt)=(0.138024,0.627629,0.89729); rgb(24pt)=(0.124814,0.645929,0.888343); rgb(25pt)=(0.111252,0.6635,0.876314); rgb(26pt)=(0.0952095,0.679829,0.859781); rgb(27pt)=(0.0688714,0.694771,0.839357); rgb(28pt)=(0.0296667,0.708167,0.816333); rgb(29pt)=(0.00357143,0.720267,0.7917); rgb(30pt)=(0.00665714,0.731214,0.766014); rgb(31pt)=(0.0433286,0.741095,0.73941); rgb(32pt)=(0.0963952,0.75,0.712038); rgb(33pt)=(0.140771,0.7584,0.684157); rgb(34pt)=(0.1717,0.766962,0.655443); rgb(35pt)=(0.193767,0.775767,0.6251); rgb(36pt)=(0.216086,0.7843,0.5923); rgb(37pt)=(0.246957,0.791795,0.556743); rgb(38pt)=(0.290614,0.79729,0.518829); rgb(39pt)=(0.340643,0.8008,0.478857); rgb(40pt)=(0.3909,0.802871,0.435448); rgb(41pt)=(0.445629,0.802419,0.390919); rgb(42pt)=(0.5044,0.7993,0.348); rgb(43pt)=(0.561562,0.794233,0.304481); rgb(44pt)=(0.617395,0.787619,0.261238); rgb(45pt)=(0.671986,0.779271,0.2227); rgb(46pt)=(0.7242,0.769843,0.191029); rgb(47pt)=(0.773833,0.759805,0.16461); rgb(48pt)=(0.820314,0.749814,0.153529); rgb(49pt)=(0.863433,0.7406,0.159633); rgb(50pt)=(0.903543,0.733029,0.177414); rgb(51pt)=(0.939257,0.728786,0.209957); rgb(52pt)=(0.972757,0.729771,0.239443); rgb(53pt)=(0.995648,0.743371,0.237148); rgb(54pt)=(0.996986,0.765857,0.219943); rgb(55pt)=(0.995205,0.789252,0.202762); rgb(56pt)=(0.9892,0.813567,0.188533); rgb(57pt)=(0.978629,0.838629,0.176557); rgb(58pt)=(0.967648,0.8639,0.16429); rgb(59pt)=(0.96101,0.889019,0.153676); rgb(60pt)=(0.959671,0.913457,0.142257); rgb(61pt)=(0.962795,0.937338,0.12651); rgb(62pt)=(0.969114,0.960629,0.106362); rgb(63pt)=(0.9769,0.9839,0.0805)}, mesh/rows=21]
table[point meta=\thisrow{c}] {%
	plots/+data/8_x_1-1.tsv};
\addlegendentry{data1}

\addplot3 [color=red, line width=2.0pt]
table[] {plots/+data/8_x_1-2.tsv};
\addlegendentry{data2}

\end{axis}
\end{tikzpicture}%
	}%
	\\
	\pgfplotsset{
		every axis post/.append style={
			width=0.4\linewidth,
			legend style={
				font=\footnotesize,
				at = {(1,0)},
				anchor = south east,
				outer sep=2pt,
			},
			every axis title/.append style={yshift=-2mm},
			every axis x label/.append style={yshift=2mm},
			yticklabel style={inner sep=0pt,xshift=-3pt},
			xticklabel style={inner xsep = 0pt},
		},%
		every axis plot/.append style={
			line width=1.5pt
		}
	}%
	{%
		\pgfplotsset{%
			every axis post/.append style={%
				execute at end axis={
					\legend{}
				},
				xmax = 0.6,
				ymin = -0.1,
				ymax = 3,
			}%
		}%
%
%
\begin{tikzpicture}

\begin{axis}[%
width=0.333\linewidth,
height=0.35\linewidth,
at={(0\linewidth,0\linewidth)},
scale only axis,
xmin=0,
xmax=1,
xlabel style={font=\color{white!15!black}},
xlabel={$t$},
ymin=0,
ymax=6,
axis background/.style={fill=white},
title style={font=\bfseries},
title={$u_0(t)$},
legend style={legend cell align=left, align=left, draw=white!15!black},
every axis title/.append style={font=\footnotesize},every axis label/.append style={font=\footnotesize},every tick label/.append style={font=\scriptsize},every x tick label/.style={},every y tick label/.style={},every z tick label/.style={},yticklabel style={/pgf/number format/fixed,/pgf/number format/precision=2},
]
\addplot [color=green]
table[]{plots/+data/4_Inputs0-1.tsv};
\addlegendentry{data1}

\addplot [color=green, dashed]
table[]{plots/+data/4_Inputs0-2.tsv};
\addlegendentry{data2}

\addplot [color=blue]
table[]{plots/+data/4_Inputs0-3.tsv};
\addlegendentry{data3}

\addplot [color=blue, dashed]
table[]{plots/+data/4_Inputs0-4.tsv};
\addlegendentry{data4}

\end{axis}
\end{tikzpicture}%
\hfill%
	}%
	{%
		\pgfplotsset{%
			every axis post/.append style={%
				xmax = 0.6,
				ymax = 0.1,
			}%
		}%
%
%
\begin{tikzpicture}

\begin{axis}[%
width=0.333\linewidth,
height=0.35\linewidth,
at={(0\linewidth,0\linewidth)},
scale only axis,
xmin=0,
xmax=1,
xlabel style={font=\color{white!15!black}},
xlabel={$t$},
ymin=-8,
ymax=0,
axis background/.style={fill=white},
title style={font=\bfseries},
title={$u_1(t)$},
legend style={legend cell align=left, align=left, draw=white!15!black},
every axis title/.append style={font=\footnotesize},every axis label/.append style={font=\footnotesize},every tick label/.append style={font=\scriptsize},every x tick label/.style={},every y tick label/.style={},every z tick label/.style={},yticklabel style={/pgf/number format/fixed,/pgf/number format/precision=2},
]
\addplot [color=green]
table[]{plots/+data/5_Inputs1-1.tsv};
\addlegendentry{$\yn=0.16$}

\addplot [color=green, dashed, forget plot]
table[]{plots/+data/5_Inputs1-2.tsv};
\addplot [color=blue]
table[]{plots/+data/5_Inputs1-3.tsv};
\addlegendentry{$\yn=0.66$}

\addplot [color=blue, dashed, forget plot]
table[]{plots/+data/5_Inputs1-4.tsv};

\end{axis}
\end{tikzpicture}%
%
	}%
	
	\caption{Profiles of the first state $w_1(\y,t)$ of the closed-loop system \eqref{eq:origSys} with 
		\eqref{eq:controlLaw} for the choices $\yn = 0.16$ and $\yn = 0.66$, which are marked by the red lines and the corresponding control efforts $u_0(t)$ and $u_1(t)$ for $\yn = 0.16$ (green) and $\yn = 0.66$ (blue). The solid lines correspond to the first element of the input $u_i(t)$, $i=0,1$, and the dashed lines to their second element.}
	\label{fig:foldingPointEffect}
\end{figure}
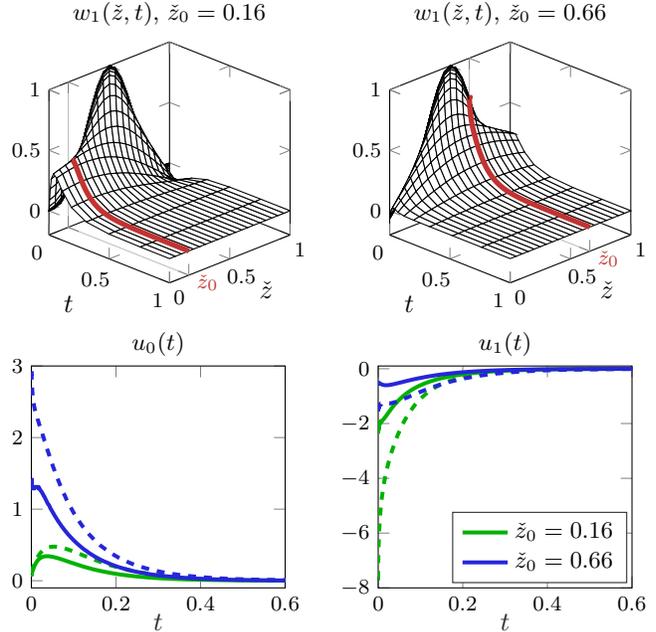

The upper plots in Figure \ref{fig:foldingPointEffect} show the resulting closed-loop profiles $w_1(\y,t)$ for the cases $\yn = 0.16$ and $\yn = 0.66$. Note that $\yn = 0.66$ requires a reordering of the states to ensure $\lambda_1 > \ldots > \lambda_{4}$. Obviously, the choice of the folding point $\yn$ has a serious influence on the spatial and temporal evolution of the state. A similar simulation result is obtained for the corresponding second state $w_2(\y,t)$. 
Moreover, 
the design parameter $\yn$ can be used to adjust the distribution of the control effort between the inputs $u_0(t)$ and $u_1(t)$. 
This can be verified in the lower plots of Figure \ref{fig:foldingPointEffect}, where it is compared in the cases $\yn = 0.16$ and $\yn = 0.66$. It can be seen that shifting the folding point to the left increases the control effort for the right input and vice versa. This effect is further investigated in \cite{Che19b}.

%
\section{Concluding Remarks}
An obvious extension in the case of coupled PDEs is to specify individual folding points for each state. This, however, leads to couplings between different spatial areas of different states in the folded system, hindering the application of the standard backstepping transformation and leading to a very involved problem.
As already mentioned, the case of non-sorted diffusion coefficients, where some or all of them are equal can be included straightforwardly. Moreover, the observer design for coupled parabolic PDEs with in-domain measurements like in \cite{Che19b} can be realized by making use of the ideas in \cite{Deu18}. The extension to time-varying reaction coefficients with the results in \cite{Ker19} poses no new challenges\cross{ for the design}. 
Furthermore, the results of the paper allow to investigate more sophisticated systems like parabolic PDEs with couplings to an ODE at an in-domain point.

Combining the folding technique with the methods for output regulation according to \cite{Deu18} and the internal model principle \cite{Deu19} to achieve robust output regulation in the bilateral case are of special interest for future work.
	Moreover, the bilateral backstepping approach opens the door to fault tolerant designs. For example, a bilateral setup 
	can be rendered fault tolerant by switching to a stabilizing unilateral controller.
	It would also be interesting to investigate the input to state stability properties of bilateral designs \wrt\ boundary disturbances.

\appendix

\bibliographystyle{plain}        
\bibliography{mybib}

\end{document}